\documentclass[titlepage,draft,12pt]{article} 
\usepackage{amssymb,amsthm} 
\usepackage[a4paper]{geometry}
\usepackage[usenames,dvipsnames]{xcolor}

\date{}

\newcommand{\ep}{\varepsilon}
\renewcommand{\qed}{{\penalty 10000\mbox{$\quad\Box$}}}
\newcommand{\re}{\mathbb{R}}

\newcommand{\n}{\mathbb{N}}

\newtheorem{thm}{Theorem}[section]
\newtheorem{thmbibl}{Theorem}

\newtheorem{rmk}[thm]{Remark}
\newtheorem{prop}[thm]{Proposition}

\newtheorem{lemma}[thm]{Lemma}
\newtheorem*{open}{Open problem}

\title{A description of all possible decay rates for solutions of 
some semilinear parabolic equations}

\author{Marina Ghisi\vspace{1ex}\\ 
{\normalsize Universit\`a degli Studi di Pisa} \\
{\normalsize Dipartimento di Matematica}\\ 
{\normalsize PISA (Italy)}\\
{\normalsize e-mail: \texttt{ghisi@dm.unipi.it}}
\and
Massimo Gobbino\vspace{1ex}\\ 
{\normalsize Universit\`a degli Studi di Pisa} \\
{\normalsize Dipartimento di Matematica}\\ 
{\normalsize PISA (Italy)}\\  
{\normalsize e-mail: \texttt{m.gobbino@dma.unipi.it}}
\and
Alain Haraux\vspace{1ex}\\ 
{\normalsize Universit\'{e} Pierre et Marie Curie} \\
{\normalsize Laboratoire Jacques-Louis Lions}\\ 
{\normalsize PARIS (France)}\\  
{\normalsize e-mail: \texttt{haraux@ann.jussieu.fr}}}

\begin{document}
\maketitle
\begin{abstract}
	We consider an abstract first order evolution equation in a
	Hilbert space in which the linear part is represented by a
	self-adjoint nonnegative operator $A$ with discrete spectrum, and
	the nonlinear term has order greater than one at the origin.  We
	investigate the asymptotic behavior of solutions.
	
	We prove that two different regimes coexist.  Close to the kernel
	of $A$ the dynamic is governed by the nonlinear term, and
	solutions (when they decay to 0) decay as negative powers of $t$.
	Close to the range of $A$, the nonlinear term is negligible, and
	solutions behave as solutions of the linearized problem.  This
	means that they decay exponentially to 0, with a rate and an
	asymptotic profile given by a simple mode, namely a one-frequency
	solution of the linearized equation.
	
	The abstract results apply to semilinear parabolic equations.
	
\vspace{6ex}

\noindent{\bf Mathematics Subject Classification 2010 (MSC2010):}
35K58, 35K90, 35B40.


\vspace{6ex}

\end{abstract}

 
\section{Introduction}

In this paper we study the asymptotic behavior of decaying solutions
to the first order evolution equation
\begin{equation}
	u'(t)+Au(t)=f(u(t))
	\quad\quad
	\forall t\geq 0,
	\label{pbm:main}
\end{equation}
where $A$ is a self-adjoint linear operator on a Hilbert space $H$ 
and $f$ is a nonlinear term.

We assume that $A$ is non-negative, but not necessarily strictly 
positive, and that its spectrum is a finite set or an increasing 
sequence of eigenvalues. We also assume that the nonlinear term has 
order greater than one in the origin, in the sense that it satisfies 
inequalities such as
\begin{equation}
	|f(u)|\leq K_{0}\left(|u|^{1+p}+|A^{1/2}u|^{1+q}\right)
	\label{pbm:order}
\end{equation}
for some positive exponents $p$ and $q$.

As a model example, we have in mind semilinear parabolic equations
such as
$$u_{t}-\Delta u+|u|^{p}u=0$$
with Neumann boundary conditions in a bounded domain 
$\Omega\subseteq\re^{n}$, or
$$u_{t}-\Delta u-\lambda_{1}(\Omega)u+|u|^{p}u=0$$
with Dirichlet boundary conditions in a bounded domain
$\Omega\subseteq\re^{n}$, where $\lambda_{1}(\Omega)$ denotes the
first eigenvalue of the Dirichlet Laplacian.  We point out that in
both cases the operator associated to the linear part has a nontrivial
kernel.

These two model examples have been investigated in the last decade in
a series of papers by the third author and collaborators.  Let us
outline their achievements and the main problems which remained
open in order to motivate our analysis.

In the Neumann case a simple application of the maximum principle
gives that all solutions decay at least as $t^{-1/p}$, even in the
$L^{\infty}$ norm.  Moreover, it is easy to exhibit examples both of
solutions decaying exactly as $t^{-1/p}$, and of solutions decaying
exponentially.  To this end, in the first case it is enough to
consider spatially homogeneous solutions $u(t,x):=v(t)$, where
necessarily $v(t)$ satisfies $v'(t)+|v(t)|^{p}v(t)=0$, and observe
that all nontrivial solutions of this ordinary differential equation
decay to 0 exactly as $t^{-1/p}$.  In the second case it is enough to
consider a symmetric domain $\Omega$, for example the interval
$(-1,1)$ in one dimension, and restrict ourselves to odd functions.
The Neumann Laplacian, when restricted to odd functions, becomes a
coercive operator, hence it is not difficult to see that all
non-trivial solutions with odd initial datum decay to 0 exponentially.

This means that in the Neumann case there is coexistence between
\emph{slow solutions}, decaying as negative powers of $t$, and
\emph{fast solutions} decaying exponentially.  It was later shown
in~\cite{ben-arbi} that actually one has the so called \emph{slow-fast
alternative}, meaning that all non-zero solutions decay to 0 either
exponentially, or exactly as $t^{-1/p}$.  Finally, in~\cite{bah-sing}
it was shown that all initial data which are small enough and close
enough to constant functions in the norm of $L^{\infty}(\Omega)$ give
rise to slow solutions.  Despite of these achievements, several
problems remained open, for example describing all possible
exponential decay rates, or even proving the existence of fast
solutions when $\Omega$ has no special symmetry, or providing an
explicit description of an open set of initial data in the norm of the
phase space $L^{2}(\Omega)$ originating slow solutions.

The Dirichlet case proved to be more difficult to tackle.  It is still
true, but more delicate to establish, that all solutions decay to 0 at
least as $t^{-1/p}$.  Moreover, a comparison argument with suitable
sub-solutions was enough to prove that all solutions with nonnegative
initial data are actually slow solutions.  Both results have been
proved in~\cite{hjk} (see also~\cite{h:nodea}).  Later on, a weak form
of slow-fast alternative was established in~\cite{bah-d}, meaning that
all non-zero solutions decay either as $t^{-1/p}$, or faster than all
negative powers of $t$.  In the same paper also the existence of
exponentially fast solutions was shown, but assuming the domain
$\Omega$ to be symmetric, and the existence of intermediate decay
rates was not excluded.  Several problems remained open, among them
existence of an open set of slow solutions, description of fast
solutions, existence of fast solutions in domains without symmetries,
and the true slow-fast alternative.

We stress that all these results were proved exploiting special
symmetries, or the existence of $A$-invariant subspaces of $H$ which
are invariant also for the nonlinear term, or comparison arguments.

What was missing is a unifying abstract theory. Filling this gap is 
the aim of this paper. So we consider the abstract evolution 
equation~(\ref{pbm:main}) and we address the following issues.
\begin{enumerate}
	\renewcommand{\labelenumi}{(\arabic{enumi})}
	\item\label{ex-slow} Existence of an open set of slow solutions.

	\item\label{ex-fast} Existence of fast solutions.

	\item\label{slow-fast} Slow-fast alternative.

	\item\label{classification} Classification of all possible decay rates.

	\item\label{description} Description of the set of solutions with
	a prescribed decay rate.
\end{enumerate}

In Theorem~\ref{thm:main-alternative} we focus on the slow-fast 
alternative. Instead of (\ref{pbm:main}) and (\ref{pbm:order}), we 
consider more generally an evolution inequality of the form
\begin{equation}
	|u'(t)+Au(t)|\leq K_{0}\left(|u(t)|^{1+p}+|A^{1/2}u(t)|^{1+q}\right)
	\quad\quad 
	\forall t\geq 0,
	\label{pbm:diff-ineq}
\end{equation}
and we prove that all its non-zero solutions, when they decay to 0,
decay either exponentially or at most as $t^{-1/p}$.  We prove also
that slow solutions move closer and closer to the kernel of $A$, in
the sense that $|A^{1/2}u(t)|$ decays faster than $|u(t)|$.  This is
clearly impossible if $A$ is strictly positive, in which case slow
solutions cannot exist.  Finally, we prove that fast solutions have an
asymptotic profile of the form $u(t)\sim v_{0}e^{-\lambda t}$, where
$\lambda$ is an eigenvalue of $A$ and $v_{0}\neq 0$ is a corresponding
eigenvector.  This settles (\ref{slow-fast}) and~(\ref{classification}).

In Theorem~\ref{thm:main-slow} we prove that slow solutions of
(\ref{pbm:main}) exist whenever $\ker(A)$ is nontrivial and the
nonlinear term satisfies (\ref {pbm:order}) and an additional sign
condition allowing the existence of global solutions.  What we
actually show is the existence of an open set of initial data
generating slow solutions, and this open set is characterized by
simple explicit inequalities such as~(\ref{hp:u0}).  This
settles~(\ref{ex-slow}).

Finally, in Theorem~\ref{thm:main-exponential} we address the
existence of fast solutions.  We prove that for every eigenvalue
$\lambda$ of $A$, and every corresponding eigenvector $v_{0}\neq 0$
which is small enough, there exists a nonempty set of initial data
originating solutions whose asymptotic profile is exactly
$v_{0}e^{-\lambda t}$.  This nonempty set is parametrized by an open
set in the subspace of $H$ generated by all eigenvectors of $A$
greater than $\lambda$, in analogy with solutions of the linearized
equation (see also Remark~\ref{rmk:exp-linear}).  This settles
(\ref{ex-fast}) and (\ref{description}).

Our results are apparently new, at least in the sense that they are
not explicitly stated elsewhere, even in the special situation where
$A$ is strictly positive and $H$ is a finite dimensional space, in
which case (\ref{pbm:main}) reduces to a system of ordinary
differential equations.  That case is usually handled by means of
Lyapunov functions or linearization theorems.  Lyapunov functions lead
to a simple proof that all solutions decay exponentially.  Classical
linearization theorems (see~\cite{hartman-1,hartman-2,hartman-3})
state that in a neighborhood of the origin the dynamic induced by the
nonlinear system is homeomorphic to the dynamic induced by the
linearized one.  Nevertheless, both methods do not lead to a
classification of decay rates, the first one because it only provides
an estimate from above, the second one because the homeomorphisms
given by the linearization theorems are just H\"{o}lder continuous,
hence they do not preserve decay rates.

When $H$ is infinite dimensional and $A$ is strictly positive, our
results seem to be new as well, at least in full generality.  The only
related literature we are aware of is~\cite{fs:nl}, where an analogous
classification of exponential decay rates has been provided for
solutions of the Navier-Stokes equation in a bounded domain, in which
case the operator $A$ is coercive and the nonlinear term is quadratic
at the origin.  As in~\cite{fs:nl}, our proof requires a careful
analysis of the asymptotic behavior of the Dirichlet
quotient~(\ref{defn:Q}).  Apart from this, our approach is quite
different, especially in the construction of fast solutions with a
prescribed asymptotic profile.  In~\cite{fs:nl} the set of such
solutions, called nonlinear spectral manifold, is characterized as a
level set of a suitable function.  That approach seems to assume that
global solutions already exist and generate a semigroup with some
regularity, assumptions which are well suited for the Navier-Stokes
equation but not for our general framework.  Therefore, what we do is
proving a stand-alone existence result through a Banach fixed point
argument.  Since fast solutions with a prescribed profile are
non-unique, it might seem impossible to obtain them by a
contraction argument.  Nevertheless, the trick is to produce them one
by one, by looking for them carefully in suitable classes of functions
where one and only one such solution is supposed to be.

In any case, the strength of our results lies mainly in dealing with
the case where the kernel of $A$ is nontrivial, which causes the
coexistence of slow and fast solutions.  This case was completely
open, even in finite dimension, apart from the Neumann and Dirichlet
examples quoted before.  When we apply our abstract results to those
examples, we obtain new proofs of the previous known results, and we
solve all open problems, at least in dimension one and two, or under a
smallness condition on $p$ depending on the dimension (this
restriction comes from the need of Sobolev embeddings in order to
verify the assumptions of our abstract results).  Actually some
results persist for any $p$, for instance the slow-fast alternative
for the Dirichlet case.  Moreover applications are not limited to the
model examples, but larger classes of semilinear parabolic equations
fit in our general framework.

From the technical point of view, both the existence of slow solutions
and the slow-fast alternative require now a careful asymptotic
analysis of what we call generalized Dirichlet quotients, defined
in~(\ref{defn:Qd}).  Similar quotients have been used also
in~\cite{ghisi:JDE2006,gg:k-decay,GGH:sol-lentes} in different
contexts (quasilinear and semilinear dissipative hyperbolic
equations), but always with the aim of estimating decay rates from
below.

We are quite optimistic about the possibility to extend our techniques
to cases where our results, as stated here, do not apply immediately.
In particular we have in mind both the model examples in any dimension
with any $p$ (for the complete result), and hyperbolic equations with
damping terms under growth conditions on the non-linearity.  These are
likely to be the directions of future investigations.

This paper is organized as follows.  In Section~\ref{sec:statements}
we clarify the functional setting, we recall two classical local
existence theorems, and we state our main abstract results, which we
prove in Section~\ref{sec:proofs}.  In Section~\ref{sec:applications}
we apply the abstract theory to semilinear parabolic problems.  

\setcounter{equation}{0}
\section{Statements}\label{sec:statements}

\subsection{Notation and classical existence results}

Throughout this paper $H$ denotes a Hilbert space, $|x|$ denotes the
norm of an element $x\in H$, and $\langle x,y\rangle$ denotes the
scalar product of two elements $x$ and $y$ in $H$.  We consider a
self-adjoint linear operator $A$ on $H$ with dense domain $D(A)$.  We
assume that $A$ is nonnegative, namely $\langle Au, u\rangle\geq 0$
for every $u \in D(A)$, so that for every $\alpha\geq 0$ the power
$A^{\alpha}u$ is defined provided that $u$ lies in a suitable domain
$D(A^{\alpha})$, which is itself a Hilbert space with norm
$$|u|_{D(A^{\alpha})} := \left(|u|^{2} +
|A^{\alpha}u|^{2}\right)^{1/2}.$$

Before stating our results, let us spend a few words on the notion of 
solution. Let us start with the linear equation
\begin{equation}
	u'(t)+Au(t)=g(t),
	\label{pbm:linear}
\end{equation}
with initial condition
\begin{equation}
	u(0)=u_{0}.
	\label{pbm:data}
\end{equation}

For our purposes we can limit ourselves to consider \emph{strong
solutions} of (\ref{pbm:linear}), namely functions $u$ defined in some
time-interval $[0,T]$ and such that for almost every $t\in(0,T)$ one 
has that $u'(t)$ exists, $u(t)\in D(A)$, and (\ref{pbm:linear}) is 
satisfied.

We recall the following classical result (we refer for example
Theorem~3.6 in~\cite{brezis} where the same regularity is obtained in
a more general nonlinear setting).

\begin{thmbibl}[Linear equation -- Existence]\label{thmbibl:linear}
	Let $H$ be a Hilbert space, and let $A$ be a self-adjoint
	nonnegative operator on $H$ with dense domain $D(A)$.  Let us
	assume that $T>0$, $g\in L^{2}((0,T),H)$ and $u_{0}\in
	D(A^{1/2})$.
	
	Then problem (\ref{pbm:linear})--(\ref{pbm:data}) has a unique 
	solution with the following regularity
	\begin{eqnarray}
		 & u\in C^{0}\left([0,T],D(A^{1/2})\right), & 
		\label{th:reg-cont}  \\
		\noalign{\vspace{1ex}}
		 & u\in W^{1,2}\left((0,T),H)\right)\cap
		L^{2}\left((0,T),D(A)\right), & 
		\label{th:reg-w12}  \\
		\noalign{\vspace{1ex}}
		 & \mbox{the function $t\to|A^{1/2}u(t)|^{2}$ is absolutely 
		continuous in $[0,T]$.} & 
		\label{th:reg-ac}
	\end{eqnarray}
\end{thmbibl}

We point out that the solution is defined as long as the forcing term
$g(t)$ is defined.  If $g\in L^{2}\left((0,T),H\right)$ for every
$T>0$, then the solution is defined for every $t\geq 0$.

We also mention that the solution provided by
Theorem~\ref{thmbibl:linear} can be represented by the well-known integral 
formula
$$u(t)=e^{-tA}u_{0}+\int_{0}^{t}e^{-(t-s)A}g(s)\,ds.$$

Now we consider a semilinear equation of the form
\begin{equation}
	u'(t)+Au(t)=f(u(t)),
	\label{pbm:nl}
\end{equation}
for which we have the following local existence result.

	\begin{thmbibl}[Semi-linear equation -- Local
	existence]\label{thmbibl:nl}
	
	Let $H$ be a Hilbert space, and let
	$A$ be a self-adjoint nonnegative operator on $H$ with dense
	domain $D(A)$.  Let $R>0$, let $B_{R}:=\left\{u\in
	D(A^{1/2}):|u|_{D(A^{1/2})}< R\right\}$, and let $f:B_{R}\to H$
	be a function.
	
	Let us assume that there exists a constant $L$ such that
	\begin{equation}
		|f(u)-f(v)|\leq L|u-v|_{D(A^{1/2})}
		\quad\quad
		\forall (u,v)\in[B_{R}]^{2}.
		\label{hp:f-lip}
	\end{equation}
	
	Then for every $u_{0}\in B_{R}$ there exist $T>0$, and a unique
	local solution $u$ to problem (\ref{pbm:nl})--(\ref{pbm:data})
	satisfying (\ref{th:reg-cont}) through (\ref{th:reg-ac}).  This
	solution can be continued to a solution defined in a maximal
	interval $[0,T_{*})$, where either $T_{*}=+\infty$ or
	\begin{equation}
		\lim_{t\to T_{*}^{-}}|u(t)|_{D(A^{1/2})}=R.
		\label{th:alternative}
	\end{equation}
	
\end{thmbibl}

The proof of Theorem~\ref{thmbibl:nl} is completely standard.  It is
enough to consider the map that associates to every $v\in
C^{0}([0,T],D(A^{1/2}))$ the solution of the linear problem
(\ref{pbm:linear})--(\ref{pbm:data}) with $g(t):=f(v(t))$.  If $T>0$
is small enough, this map turns out to be a contraction, and the
unique fixed point is the required (strong) local solution.

\subsection{Main results}

Our first result provides a classification of all possible decay 
rates for decaying solutions to differential inequalities of the 
form~(\ref{pbm:diff-ineq}).

\begin{thm}[Classification of decay rates]\label{thm:main-alternative}
	Let $H$ be a Hilbert space, and let $A$ be a self-adjoint nonnegative
	operator on $H$ with dense domain $D(A)$. 
	
	Let $g\in L^{2}\left((0,T),H\right)$ for every $T>0$, and let
	$u\in C^{0}\left([0,+\infty),D(A^{1/2})\right)$ be a global
	solution of (\ref{pbm:linear}) in the sense of
	Theorem~\ref{thmbibl:linear}.  Let us assume that
	\begin{enumerate}
		\renewcommand{\labelenumi}{(\roman{enumi})} 
		
		\item the spectrum of $A$ is a finite set or an increasing
		sequence of eigenvalues,
	
		\item  $u$ is a decaying solution in the sense that
		\begin{equation}
			\lim_{t\to +\infty}|u(t)|_{D(A^{1/2})}=0,
			\label{hp:u-limit}
		\end{equation}
	
		\item there exist $p>0$, $q>0$, and $K_{0}\geq 0$ such that
		\begin{equation}
			\left|g(t)\right|\leq K_{0}
			\left(|u(t)|^{1+p}+|A^{1/2}u(t)|^{1+q}\right)
			\quad\quad
			\forall t\geq 0.
			\label{hp:diff-ineq}
		\end{equation}
	\end{enumerate}
		
	Then one and only one of the following statements apply.
	\begin{enumerate}
		\renewcommand{\labelenumi}{(\arabic{enumi})}
		
		\item \emph{(Null solution)} The solution is the zero-solution
		$u(t)\equiv 0$ for every $t\geq 0$.
		
		\item \emph{(Slow solutions)} There exist positive constants 
		$M_{1}$ and $M_{2}$ such that 
		\begin{equation}
			|u(t)|\geq \frac{M_{1}}{(1+t)^{1/p}}
			\quad\quad
			\forall t\geq 0,
			\label{th:slow}
		\end{equation}
		\begin{equation}
			|A^{1/2}u(t)|\leq M_{2}|u(t)|^{1+p}
			\quad\quad
			\forall t\geq 0.
			\label{th:range}
		\end{equation}
		
		\item \emph{(Spectral fast solutions)} There exist an
		eigenvalue $\lambda>0$ of $A$, and a corresponding eigenvector
		$v_{0}\neq 0$, such that
		\begin{equation}
			\lim_{t\to+\infty}\left|u(t)-v_{0}e^{-\lambda t}\right|_{D(A^{1/2})}
			e^{\gamma t}=0
			\label{th:fast+}
		\end{equation}
		for every 
		\begin{equation}
			\gamma<\min\left\{\strut\beta,(1+p)\lambda,
			(1+q)\lambda\right\},
			\label{defn:delta}
		\end{equation}
		where $\beta=+\infty$ if the spectrum of $A$ is finite and
		$\lambda$ is its maximum, and $\beta$ is the smallest
		eigenvalue of $A$ larger than $\lambda$ otherwise.
	\end{enumerate}
\end{thm}

\begin{rmk}
	\begin{em}
		When the kernel of $A$ is non-trivial, a differential
		inequality such as~(\ref{pbm:diff-ineq}) does not guarantee
		that all its solutions in a neighborhood of the origin tend to
		0 (just think to the ordinary differential equation
		$u'=u^{3}$).  This is the reason why we need
		assumption~(\ref{hp:u-limit}).
		
		In other words, there might be coexistence of solutions that
		decay to 0 and solutions that do not decay, or even do not
		globally exist.  When this is the case, our result classifies
		all possible decay rates of decaying solutions, regardless of
		non-decaying ones.
	\end{em}
\end{rmk}

\begin{rmk}
	\begin{em}
		Concerning the null solution,
		Theorem~\ref{thm:main-alternative} extends the well-known
		backward uniqueness results of the seminal 
		papers~\cite{bt,ghidaglia}.
		Assuming forward uniqueness, which holds true for large
		classes of equations, classical backward uniqueness results
		read as follows.  If $u(t)=0$ for some $t\geq 0$ (hence also
		for all subsequent times), then $u(t)=0$ for all $t\geq 0$.
		Our result extends the classical one by showing that, if
		$u(t)$ decays at infinity faster than $e^{-ct}$ for all $c>0$,
		then $u(t)=0$ for all $t\geq 0$.
	\end{em}
\end{rmk}

\begin{rmk}
	\begin{em}
		Concerning slow solutions, we point out that only the exponent
		$p$ in (\ref{hp:diff-ineq}) appears in the decay rate, while
		$q$ is irrelevant provided it is positive.  Roughly speaking,
		this happens because slow solutions move closer and closer to
		the kernel of $A$, as suggested by the otherwise unnatural
		estimate (\ref{th:range}) in which $|A^{1/2}u(t)|$ is
		controlled with a higher power of $|u(t)|$.  Close to the
		kernel of $A$, the term $|A^{1/2}u(t)|$ can be neglected, and
		this justifies the disappearance of $q$ in the final decay
		rate.  
		
		Let us write $u(t)$ as the sum of its projection $P_{K}u(t)$ 
		into $\ker(A)$, and its ``range component'' $u(t)-P_{K}u(t)$ 
		orthogonal to $\ker(A)$. Since the operator $A$ is coercive 
		when restricted to the range of $A$, estimate 
		(\ref{th:range}) implies that there exists a constant $c$ 
		such that
		$$\left|u(t)-P_{K}u(t)\right|\leq
		c|A^{1/2}u(t)|\leq cM_{2}|u(t)|^{1+p}.$$
		
		In other words, when $u(t)$ decays to 0, its range component
		always decays faster, so that the slowness of $u(t)$ is due
		uniquely to its component $P_{K}u(t)$ with respect to the
		kernel.  This is consistent with previous results
		(see~\cite{h:nodea}).  This shows also that slow solutions
		cannot exist when the operator $A$ is coercive.
		
		The exponent $(1+p)$ in (\ref{th:range}) is optimal.  This can
		be seen by considering the case where $H=\re^{2}$, $p=2$, and the
		evolution problem reduces to the following system of ordinary
		differential equations 
		$$\left\{
		\begin{array}{l}
			x'(t)= -x^{3}(t),  \\
			y'(t)+ y(t)=x^{3}(t).
		\end{array}
		\right.$$
		
		A solution of the first equation is $x(t)=(1+2t)^{-1/2}$.  At
		this point it is possible to prove that all solutions of the
		second equation decay as the forcing term, hence as 
		$(1+2t)^{-3/2}$. Therefore, in this example we have that 
		$|u(t)|\sim|x(t)|\sim(1+2t)^{-1/2}$, while 
		$|A^{1/2}u(t)|=|y(t)|\sim(1+2t)^{-3/2}=|u(t)|^{1+p}$.

	\end{em}
\end{rmk}

\begin{rmk}
	\begin{em}
		Concerning fast solutions, the possible asymptotic profiles 
		are described by~(\ref{th:fast+}), which also provides an 
		estimate for the remainder. We point out that 
		(\ref{defn:delta}) is optimal. This can be seen by 
		considering the case where $H=\re^{2}$ and the evolution 
		problem reduces to the following system of ordinary 
		differential equations
		$$\left\{
		\begin{array}{l}
			x'(t)+\lambda x(t)=0,  \\
			y'(t)+\beta y(t)=|x(t)|^{1+p}+|x(t)|^{1+q}.
		\end{array}
		\right.$$
		
		A solution of the first equation is $x(t)=e^{-\lambda t}$. At 
		this point, solutions of the second equation can decay as 
		$e^{-\eta t}$, where $\eta$ is the right-hand side 
		of~(\ref{defn:delta}), or even as $te^{-\eta t}$ in case of 
		resonance. 
	\end{em}
\end{rmk}

We conclude this long discussion on Theorem~\ref{thm:main-alternative} by 
mentioning the following.

\begin{open}
	\begin{em}
		Is it possible to weaken assumption (\ref{hp:u-limit}) by
		asking just that $|u(t)|\to 0$ as $t\to +\infty$, namely by
		requiring the limit in $H$ instead of $D(A^{1/2})$?  Our
		proof requires the assumption as stated, but we have no
		counterexamples with the weaker requirement.  Actually we have
		no examples at all of solutions that decay to 0 in $H$ but not
		in $D(A^{1/2})$.
	\end{em}
\end{open}

The slow-fast alternative alone does not guarantee the existence of
both slow solutions and fast solutions.  Next result provides
sufficient conditions for the existence of an open set of slow
solutions.

\begin{thm}[Existence of slow solutions]\label{thm:main-slow}
	
	Let $H$ be a Hilbert space, and let $A$ be a self-adjoint nonnegative
	operator on $H$ with dense domain $D(A)$. Let $f:B_{R}\to H$ be a 
	function, with $R>0$ and $B_{R}$ as in Theorem~\ref{thmbibl:nl}.
	
	Let us assume that
	\begin{enumerate}
		\renewcommand{\labelenumi}{(\roman{enumi})}
		
		\item $\ker(A)\neq \{0\}$, and there exists a constant $\nu>0$
		such that 
		\begin{equation}
			|Au|^{2}\geq\nu|A^{1/2}u|^{2}
			\quad\quad
			\forall u\in D(A),
			\label{hp:A-nu}
		\end{equation}
		
		\item there exists a constant $L$ such that (\ref{hp:f-lip})
		holds true,
		
		\item there exist $p>0$, $q>0$, and $K_{0}\geq 0$ such that
		\begin{equation}
			|f(u)|\leq K_{0}\left(|u|^{1+p}+|A^{1/2}u|^{1+q}\right)
			\quad\quad
			\forall u\in B_{R},
			\label{hp:f-ho}
		\end{equation}
		and in addition
		\begin{equation}
			\langle u,f(u)\rangle\leq 0
			\quad\quad
			\forall u\in B_{R}.
			\label{hp:f-sign}
		\end{equation}
\end{enumerate}

	Then there exists an open set $\mathcal{S}\subseteq B_{R}$ (open
	with respect to the norm of $D(A^{1/2})$) with the following
	property.  For every $u_{0}\in\mathcal{S}$, the unique solution
	$u$ of (\ref{pbm:nl})--(\ref{pbm:data}) provided by
	Theorem~\ref{thmbibl:nl} is actually global, and slow in the sense
	that it satisfies (\ref{th:slow}).
\end{thm}

\begin{rmk}
	\begin{em}
		Let us briefly comment on the hypotheses of
		Theorem~\ref{thm:main-slow}.  Concerning the operator $A$, we
		already pointed out that $\ker(A)\neq \{0\}$ is a necessary
		condition for the existence of slow solutions, while
		(\ref{hp:A-nu}) is automatic if the spectrum of $A$ is a
		finite set or an increasing sequence of eigenvalues.
		
		Concerning the nonlinear term, assumption (\ref{hp:f-lip})
		comes from the local existence result, while (\ref{hp:f-ho})
		means that the nonlinear term has order higher than one at the
		origin, in accordance with (\ref{hp:diff-ineq}).
		
		Assumption~(\ref{hp:f-sign}) guarantees that the function
		$t\to|u(t)|$ is nonincreasing, and this is exploited in the
		proof in order to keep the solution inside $B_{R}$.  This
		assumption can be weakened in several ways, for example by
		requiring only that $\langle u,f(u)\rangle\leq|A^{1/2}u|^{2}$
		for every $u\in B_{R}$, but it can not be dropped completely.
		Indeed this is a sort of sign condition, and when it is
		violated one can not guarantee even the existence of global
		solutions, as in the case of the ordinary differential
		equation $u'=u^{3}$.
	\end{em}
\end{rmk}

\begin{rmk}
	\begin{em}
		Concerning the conclusion of Theorem~\ref{thm:main-slow}, we
		point out that we prove the existence of an open set of
		solutions decaying \emph{at most} as $t^{-1/p}$.  In general
		it is not true that solutions decay exactly as $t^{-1/p}$.  As
		a matter of fact, they can even not to decay at all.  For
		example, the assumptions are satisfied in the extreme case
		where both $A$ and the nonlinear term are identically~0, and
		in that case all solutions are stationary (which implies
		slow).
		
		Of course, when we know that a solution is slow, we can always
		apply Theorem~\ref{thm:main-alternative} and deduce that it
		satisfies (\ref{th:range}) in addition to~(\ref{th:slow}).
	\end{em}
\end{rmk}

Our last result concerns the existence of families of fast solutions.

\begin{thm}[Existence of fast solutions]
	\label{thm:main-exponential}
	
	Let $H$ be a Hilbert space, and let $A$ be a self-adjoint nonnegative
	operator on $H$ with dense domain $D(A)$. Let $f:B_{R}\to H$ be a 
	function, with $R>0$ and $B_{R}$ as in Theorem~\ref{thmbibl:nl}.
	
	Let us assume that
	\begin{enumerate}
		\renewcommand{\labelenumi}{(\roman{enumi})}
		
		\item the spectrum of $A$ is a finite set or an increasing
		sequence of eigenvalues,
		
		\item there exist $p>0$ and $L\geq 0$ such that
		\begin{equation}
			|f(u)-f(v)|\leq L\left(
			|u|^{p}_{D(A^{1/2})}+|v|^{p}_{D(A^{1/2})}\right)|u-v|_{D(A^{1/2})}
			\label{hp:f-plip}
		\end{equation}
		for every $u$ and $v$ in $B_{R}$, and in addition
		\begin{equation}
			f(0)=0.
			\label{hp:f0}
		\end{equation}
	\end{enumerate}

	Let $\lambda>0$ be an eigenvalue of $A$, and let $H=H_{-}\oplus
	H_{+}$ be the orthogonal decomposition of $H$ where $H_{-}$ is the
	closure of the subspace generated by all eigenvectors of $A$
	relative to eigenvalues less than or equal to $\lambda$, and
	$H_{+}$ is the closure of the subspace generated by all
	eigenvectors of $A$ relative to eigenvalues greater than
	$\lambda$.
	
	Then there exists $r_{0}>0$ with the following property. For 
	every eigenvector $v_{0}$ relative to $\lambda$, and every 
	$w_{0}\in H_{+}\cap D(A^{1/2})$ such that
	\begin{equation}
		|v_{0}|_{D(A^{1/2})}+
		|w_{0}|_{D(A^{1/2})}\leq r_{0},
		\label{hp:smallness}
	\end{equation}
	there exists $w_{1}\in H_{-}$ such that the unique local solution,
	provided by Theorem~\ref{thmbibl:nl},
	to problem (\ref{pbm:nl})--(\ref{pbm:data}) with initial condition
	$u_{0}:=w_{0}+w_{1}$ is
	actually global and
	\begin{equation}
		\lim_{t\to +\infty}
		\left|e^{\lambda t}u(t)-v_{0}\right|_{D(A^{1/2})}=0.
		\label{th:limit}
	\end{equation}

\end{thm}

\begin{rmk}
	\begin{em}
		Assumptions (\ref{hp:f-plip}) and (\ref{hp:f0}) imply both
		(\ref{hp:f-lip}) and (\ref{hp:f-ho}) with $p=q$.  Moreover,
		(\ref{hp:f-plip}) is stronger than (\ref{hp:f-lip}) because it
		requires that the local Lipschitz constant of $f$ vanishes at
		the origin.  We emphasize that we do not impose any sign
		condition on $f$, and therefore the assumptions of
		Theorem~\ref{thm:main-exponential} are not enough to guarantee
		the existence of a global solution for every $u_{0}\in B_{R}$.
		For this reason, also the global existence part of the
		statement is nontrivial.
		
		We observe also that (\ref{hp:f-plip}) could be stated with two 
		different exponents as follows
		$$|f(u)-f(v)|\leq L\left(
		|u|^{p}+|v|^{p}+|A^{1/2}u|^{q}+|A^{1/2}v|^{q}\right)
		|u-v|_{D(A^{1/2})},$$
		or even with four exponents, but this would be useless because
		the exponents do not appear in the conclusion.  What is
		relevant here is just that they are both positive, and thus
		there is no loss of generality in assuming that they are
		equal.
		
		Of course, when we know that a solution is fast, we can always
		apply Theorem~\ref{thm:main-alternative} and deduce that the
		remainder in (\ref{th:limit}) satisfies~(\ref{th:fast+}).
	\end{em}
\end{rmk}

\begin{rmk}\label{rmk:exp-linear}
	\begin{em}
		Let $\{\lambda_{k}\}$ denote the increasing sequence of
		eigenvalues of $A$, which only for simplicity we assume of
		multiplicity one, and let $\{e_{k}\}$ denote a corresponding
		orthonormal system.  All solutions of the linear homogeneous
		equation $u'(t)+Au(t)=0$ can be represented as
		$$u(t)=\sum_{k=0}^{\infty}u_{0k}e^{-\lambda_{k}t}e_{k},$$
		where $\{u_{0k}\}$ are the components of the initial
		condition $u_{0}$ with respect to $\{e_{k}\}$.
		
		If we fix an eigenvalue $\lambda_{i}$, then the solutions
		decaying as $e^{-\lambda_{i}t}$ (up to a multiplicative
		constant) are those of the form
		$$u(t)=u_{0i}e^{-\lambda_{i}t}e_{i}+
		\sum_{k=i+1}^{\infty}u_{0k}e^{-\lambda_{k}t}e_{k},$$
		hence they are parametrized by the eigenvector
		$u_{0i}e_{i}$, which appears in the term which gives the
		asymptotic profile, and by the projection of the initial
		condition in the space $H_{+}$ generated by all eigenvectors
		relative to eigenvalues greater than $\lambda_{i}$.
		
		Theorem~\ref{thm:main-exponential} shows that the same
		parameters are involved in the nonlinear case, provided that
		we restrict to a neighborhood of the origin.  Roughly
		speaking, what we prove is an existence result for
		solutions of~(\ref{pbm:nl}) satisfying a mix of conditions at
		$t=0$ and $t=+\infty$, namely
		\begin{itemize}
			\item  a prescribed asymptotic profile (in a certain 
			sense a condition at $t=+\infty$),
		
			\item  a prescribed component of the initial datum 
			$u_{0}$ with respect to the subspace $H_{+}$.
		\end{itemize}
		
		This suggests also that solutions with decay rate exactly
		$e^{-\lambda_{i}t}$ (up to multiplicative constants) are
		``generic'' among solutions with decay rate at least
		$e^{-\lambda_{i}t}$ and, when the kernel of $A$ in
		non-trivial, slow solutions are ``generic'' among all decaying
		solutions.  This rough idea that slower behaviors are always
		``generic'', in a sense to be made precise, would probably
		deserve further investigation in the future.
	\end{em}
\end{rmk}

\setcounter{equation}{0}
\section{Proofs}\label{sec:proofs}

\subsection{Estimates for linear equations}

Let us start with two simple estimates for differential 
inequalities and integrals.

\begin{lemma}\label{lemma:z}
	Let $\varphi:[0,+\infty)\to[0,+\infty)$ be a continuous function 
	such that
	\begin{equation}
		\lim_{t\to +\infty}\varphi(t)=0.
		\label{hp:phi-to-0}
	\end{equation}
	
	Let $c>0$, and let $z:[0,+\infty)\to[0,+\infty)$ be an absolutely 
	continuous function such that
	\begin{equation}
		z'(t)\leq -cz(t)+\varphi(t)
		\label{hp:ineq-z}
	\end{equation}
	for almost every $t>0$.
	
	Then we have that
	\begin{equation}
		\lim_{t\to +\infty}z(t)=0.
		\label{th:z-to-0}
	\end{equation}
\end{lemma}

\paragraph{\textmd{\textit{Proof}}}

Integrating the differential inequality (\ref{hp:ineq-z}) it follows 
that
$$z(t)\leq z(0)e^{-ct}+e^{-ct}\int_{0}^{t}e^{cs}\varphi(s)\,ds.$$

The first term in the right-hand side tends to 0 as $t\to +\infty$.
For the second term we can reason as follows.  First $\varphi(t)$, being
continuous and convergent at infinity, is bounded by some $M>0$.  Then
for any $\varepsilon>0$ there exists $T(\varepsilon)\geq 0$ such that
$$\varphi(t)\leq \varepsilon
\quad\quad
\forall t\geq T(\varepsilon).$$ 

By splitting the integral on the two sub-intervals $[0,
T(\varepsilon)]$ and $[T(\varepsilon),t]$ we find
$$e^{-ct}\int_{0}^{t}e^{cs}\varphi(s)\,ds \leq \frac{M}{c}
e^{c(T(\varepsilon)-t)} + \frac{\varepsilon}{c}
\quad\quad
\forall t\geq T(\varepsilon),$$
so that (\ref{th:z-to-0}) follows by letting first $t\to +\infty$ and
then $\ep\to 0^{+}$.\qed

\begin{lemma}\label{lemma:int}
	Let $\delta>0$, and let $\varphi:[0,+\infty)\to[0,+\infty)$ be a
	continuous function such that
	\begin{equation}
		\lim_{t\to +\infty}\varphi(t)e^{\gamma t}=0
		\quad\quad
		\forall\gamma<\delta.
		\label{hp:phi-gamma}
	\end{equation}
	
	Then for every $\alpha<\delta$ we have that the integral
	\begin{equation}
		\int_{0}^{+\infty}e^{\alpha s}\varphi(s)\,ds
		\label{th:int-conv}
	\end{equation}
	converges, and
	\begin{equation}
		\lim_{t\to +\infty}e^{\gamma t}
		\int_{t}^{+\infty}e^{\alpha(s-t)}\varphi(s)\,ds=0
		\quad\quad
		\forall\gamma<\delta.
		\label{th:int-gamma}
	\end{equation}

\end{lemma}

\paragraph{\textmd{\textit{Proof}}}

Let us choose $\eta\in(\alpha,\delta)$.  Due to (\ref{hp:phi-gamma}),
there exists a constant $c_{\eta}$ such that $\varphi(s)\leq
c_{\eta}e^{-\eta s}$ for every $s\geq 0$, which easily implies the
convergence of~(\ref{th:int-conv}).  At this point
(\ref{th:int-gamma}) is obvious if $\gamma\leq\alpha$.  On the other
hand, for $\gamma\in (\alpha, \delta)$ we can write 
$$e^{\gamma t}\int_{t}^{+\infty}e^{\alpha(s-t)}\varphi(s)\,ds
=\int_{t}^{+\infty}e^{(\gamma -\alpha)(t-s)} e^{\gamma
s}\varphi(s)\,ds \leq \frac{1}{\gamma -\alpha} \sup_{s\geq t}
e^{\gamma s}\varphi(s)$$
for every $t\geq 0$,
so that in this case (\ref{th:int-gamma}) follows from assumption
(\ref{hp:phi-gamma}).\qed \medskip

Now we prove estimates for solutions to the 
non-homogeneous linear equation
\begin{equation}
	w'(t)+Aw(t)=\psi(t)
	\quad\quad
	\forall t\geq 0.
	\label{hp:ODE-X}
\end{equation}

Here we assume that $X$ is a Hilbert space, $A$ is a self-adjoint
nonnegative linear operator on $X$ with dense domain $D(A)$, the
forcing term $\psi$ is in $L^{2}((0,+\infty),X)$, and $w\in
C^{0}([0,+\infty),D(A^{1/2}))$ is a solution of~(\ref{hp:ODE-X})
in the sense of Theorem~\ref{thmbibl:linear}.  When we apply these
estimates in the proof of the main results, $X$ is a suitable subspace
of $H$, different from case to case, and (\ref{hp:ODE-X}) is the
projection of~(\ref{pbm:linear}) onto $X$.

\begin{lemma}[Supercritical frequences]\label{lemma:ODE-super}
	Let us assume that there exist $\beta>0$ and $\delta> 0$ 
	such that
	\begin{equation}
		|A^{1/2}x|^{2}\geq\beta|x|^{2}
		\quad\quad
		\forall x\in D(A^{1/2}),
		\label{hp:ODE-spA}
	\end{equation}
	\begin{equation}
		\lim_{t\to+\infty}|\psi(t)|e^{\gamma t}=0
	   \quad\quad
	   \forall \gamma<\delta.
	   \label{hp:ODE-sppsi}  
	\end{equation}
	
	Then we have that 
	\begin{equation}
		\lim_{t\to+\infty}|w(t)|_{D(A^{1/2})}e^{\gamma t}=0
	   \quad\quad
	   \forall \gamma<\min\{\beta,\delta\}.
		\label{th:ODE-sp}
	\end{equation}
\end{lemma}

\paragraph{\textmd{\textit{Proof}}}

Let us consider the function $F_{\gamma}(t):=e^{2\gamma
t}|A^{1/2}w(t)|^{2}$.  Due to (\ref{hp:ODE-spA}), the norm of $w(t)$
in $D(A^{1/2})$ is equivalent to $|A^{1/2}w(t)|$, hence
(\ref{th:ODE-sp}) is equivalent to proving that $F_{\gamma}(t)\to 0$
as $t\to+\infty$ for every $\gamma<\min\{\beta,\delta\}$.

To this end, we choose $\ep\in(0,\beta-\gamma)$ and we estimate the
time-derivative as follows
\begin{eqnarray*}
	F_{\gamma}'(t) & = & 
	-2e^{2\gamma t}|Aw(t)|^{2}+
	2e^{2\gamma t}\langle Aw(t),\psi(t)\rangle+
	2\gamma e^{2\gamma t}|A^{1/2}w(t)|^{2}   \\
	\noalign{\vspace{1ex}}
	 & \leq & -2e^{2\gamma t}|Aw(t)|^{2}+
	 \frac{2\ep}{\beta}e^{2\gamma t}|Aw(t)|^{2}+
	 \frac{\beta}{2\ep}e^{2\gamma t}|\psi(t)|^{2}+
	2\gamma e^{2\gamma t}|A^{1/2}w(t)|^{2}.
\end{eqnarray*}

From assumption (\ref{hp:ODE-spA}) we have that 
$|Aw(t)|^{2}\geq\beta|A^{1/2}w(t)|^{2}$, hence
$$F_{\gamma}'(t)\leq
-2(\beta-\gamma-\ep)F_{\gamma}(t)+
\frac{\beta}{2\ep}e^{2\gamma t}|\psi(t)|^{2}.$$

Now let us set
$$z(t):=F_{\gamma}(t),
\hspace{3em}
c:=2(\beta-\gamma-\ep),
\hspace{3em}
\varphi(t):=\frac{\beta}{2\ep}e^{2\gamma t}|\psi(t)|^{2}.$$

Since $\gamma<\delta$, assumption (\ref{hp:ODE-sppsi})
implies (\ref{hp:phi-to-0}).  Therefore, we can apply
Lemma~\ref{lemma:z} and deduce that $F_{\gamma}(t)\to 0$ as
$t\to+\infty$, which completes the proof.\qed

\begin{lemma}[Subcritical frequences]\label{lemma:ODE-sub}
	Let us assume that there exist $\delta>\alpha\geq 0$ such that
	\begin{equation}
		|A^{1/2}x|^{2}\leq\alpha|x|^{2}
		\quad\quad
		\forall x\in X,
		\label{hp:ODE-sbA}
	\end{equation}
	\begin{equation}
		\lim_{t\to+\infty}|\psi(t)|e^{\gamma t}=0
	   \quad\quad
	   \forall \gamma<\delta.
		\label{hp:ODE-sbpsi}
	\end{equation}

	Then the following limit
	\begin{equation}
		x_{0}:=\lim_{t\to+\infty}e^{tA}w(t)
		\label{th:defn-x0}
	\end{equation}
	exists, and
	\begin{equation}
		\lim_{t\to+\infty}
		\left|w(t)-e^{-tA}x_{0}\right|_{D(A^{1/2})}e^{\gamma t}=0
		\quad\quad
		\forall \gamma<\delta.
		\label{th:remainder}
	\end{equation}
\end{lemma}

\paragraph{\textmd{\textit{Proof}}}

Every solution of (\ref{hp:ODE-X}) is given by the explicit formula
\begin{equation}
	w(t)=e^{-tA}\left(w(0)+
	\int_{0}^{t}e^{sA}\psi(s)\,ds\right)
	\quad\quad
	\forall t\geq 0.
	\label{eq:ODE-sexpl}
\end{equation}

We claim that the integral in the right-hand side has a finite limit
when $t\to +\infty$.  Indeed assumption (\ref{hp:ODE-sbA}) guarantees
that $e^{sA}$ is a bounded operator on $X$ with norm less than or
equal to $e^{\alpha s}$, hence it is enough to prove that the integral
$$\int_{0}^{+\infty}e^{\alpha s}|\psi(s)|\,ds$$
converges.  Since $\alpha<\delta$, this follows from
Lemma~\ref{lemma:int} applied with $\varphi(t):=|\psi(t)|$, and proves
(\ref{th:defn-x0}) with
$$x_{0}:=w(0)+\int_{0}^{+\infty}e^{sA}\psi(s)\,ds.$$

Now (\ref{eq:ODE-sexpl}) can be rewritten as
$$w(t)=e^{-tA}\left(x_{0}-\int_{t}^{+\infty}e^{sA}\psi(s)\,ds\right)=
e^{-tA}x_{0}-\int_{t}^{+\infty}e^{(s-t)A}\psi(s)\,ds.$$

Exploiting again (\ref{hp:ODE-sbA}), we have now that
$$\left|w(t)-e^{-tA}x_{0}\right|e^{\gamma t}\leq
e^{\gamma t}\int_{t}^{+\infty}\left|e^{(s-t)A}\psi(s)\right|\,ds
\leq e^{\gamma t}\int_{t}^{+\infty}e^{\alpha(s-t)}
\left|\psi(s)\right|\,ds,$$
so that (\ref{th:remainder}) follows from conclusion
(\ref{th:int-gamma}) of Lemma~\ref{lemma:int}, applied once again with
$\varphi(t):=|\psi(t)|$ (we remind that the norm in $H$ and
$D(A^{1/2})$ are in this case equivalent owing to
(\ref{hp:ODE-sbA})).\qed

\subsection{Generalized Dirichlet quotients}

In this section we consider the classical \emph{Dirichlet quotient}
\begin{equation}
	Q(t):=\frac{|A^{1/2}u(t)|^{2}}{|u(t)|^{2}}.
	\label{defn:Q}
\end{equation}

We also consider the following \emph{generalized Dirichlet quotient}
\begin{equation}
	Q_{d}(t):=\frac{|A^{1/2}u(t)|^{2}}{|u(t)|^{2+d}},
	\label{defn:Qd}
\end{equation}
defined for every $d\geq 0$.

The aim of the next result is providing estimates for the 
time-derivative of $Q(t)$ and $Q_{d}(t)$ when $u$ is a solution of 
a linear equation such as~(\ref{pbm:linear}).

\begin{lemma}[Time-derivatives of Dirichlet quotients]
	Let $H$ be a Hilbert space, and let $A$ be a self-adjoint
	nonnegative operator on $H$ with dense domain $D(A)$.  Let
	$(a,b)\subseteq(0,+\infty)$ be an interval, let $g\in
	L^{2}((a,b),H)$, and let $u\in C^{0}((a,b),D(A^{1/2}))$ be a
	solution of~(\ref{pbm:linear}) in $(a,b)$ in the sense of
	Theorem~\ref{thmbibl:linear}.
	
	Let us assume that $u(t)\neq 0$ for every $t\in(a,b)$.
	
	Then we have the following conclusions.
	\begin{enumerate}
		\renewcommand{\labelenumi}{(\roman{enumi})} 
		
		\item  The Dirichlet quotient $Q(t)$ defined by 
		(\ref{defn:Q}) is absolutely continuous in $(a,b)$, and
		\begin{equation}
			Q'(t)\leq
			-\frac{\left|Au(t)-Q(t)u(t)\right|^{2}}{|u(t)|^{2}}+
			\frac{|g(t)|^{2}}{|u(t)|^{2}}
			\label{th:deriv-Q}
		\end{equation}
		for almost every $t\in(a,b)$.
	
		\item Let us assume that (\ref{hp:A-nu}) holds true for some
		constant $\nu>0$.  Then for every $d>0$ the
		generalized Dirichlet quotient $Q_{d}(t)$ defined by
		(\ref{defn:Qd}) is absolutely continuous in $(a,b)$, and
		\begin{equation}
			Q_{d}'(t)\leq 
			-\nu Q_{d}(t)+2(2+d)|u(t)|^{d}\cdot|Q_{d}(t)|^{2}+
			(3+d)\frac{|g(t)|^{2}}{|u(t)|^{2+d}}
			\label{th:deriv-Qd}
		\end{equation}
		for almost every $t\in(a,b)$.
	\end{enumerate}
\end{lemma}

\paragraph{\textmd{\textit{Proof}}}

The time-derivative of (\ref{defn:Q}) is
$$Q'(t)=-2\frac{\left|Au(t)-Q(t)u(t)\right|^{2}}{|u(t)|^{2}}+
2\frac{\langle Au(t)-Q(t)u(t),g(t)\rangle}{|u(t)|^{2}}.$$

Since
$$2\langle Au(t)-Q(t)u(t),g(t)\rangle\leq
\left|Au(t)-Q(t)u(t)\right|^{2}+|g(t)|^{2},$$
estimate (\ref{th:deriv-Q}) easily follows.

The time-derivative of (\ref{defn:Qd}) is
\begin{eqnarray*}
	Q_{d}'(t) & = & -2\frac{|Au(t)|^{2}}{|u(t)|^{2+d}}
	+(2+d)\frac{|A^{1/2}u(t))|^{2}}{|u(t)|^{2}}\cdot Q_{d}(t)  \\
	\noalign{\vspace{1ex}}
	 & & +2\frac{\langle Au(t),g(t)\rangle}{|u(t)|^{2+d}}-
	(2+d)\frac{\langle Q_{d}(t)u(t),g(t)\rangle}{|u(t)|^{2}}  \\
	\noalign{\vspace{1ex}}
	 & =: & I_{1}+I_{2}+I_{3}+I_{4}.
\end{eqnarray*}

Now it is easy to see that
$$I_{2}=(2+d)|u(t)|^{d}\cdot [Q_{d}(t)]^{2},
\hspace{4em}
I_{3}\leq\frac{|Au(t)|^{2}}{|u(t)|^{2+d}}+\frac{|g(t)|^{2}}{|u(t)|^{2+d}},$$
and
\begin{eqnarray*}
	I_{4} & \leq & \frac{2+d}{|u(t)|^{2}}\left(
	|u(t)|^{d}\cdot [Q_{d}(t)]^{2}|u(t)|^{2}+
	\frac{1}{|u(t)|^{d}}\cdot|g(t)|^{2}\right) \\
	\noalign{\vspace{1ex}}
	 & = & (2+d)|u(t)|^{d}\cdot[Q_{d}(t)]^{2}+
	(2+d)\frac{|g(t)|^{2}}{|u(t)|^{2+d}},
\end{eqnarray*}
so that
$$Q_{d}'(t)\leq -\frac{|Au(t)|^{2}}{|u(t)|^{2+d}}+
2(2+d)|u(t)|^{d}\cdot[Q_{d}(t)]^{2}+
(3+d)\frac{|g(t)|^{2}}{|u(t)|^{2+d}}.$$
	
The first term in the right-hand side is less than or equal to $-\nu
Q_{d}(t)$ because of (\ref{hp:A-nu}), and this proves
(\ref{th:deriv-Qd}).\qed

\subsection{Proof of Theorem~\ref{thm:main-alternative}}

Let us describe the scheme of the proof before entering into details.
In the first section of the proof we get rid of the null solution.
Indeed we prove that $u(T)=0$ for some $T\geq 0$ if and only if
$u(t)=0$ for every $t\geq 0$.  This is a result of forward and
backward uniqueness of the null solution.  After proving it, we can
assume that
\begin{equation}
	u(t)\neq 0
	\quad\quad
	\forall t\geq 0,
	\label{hp:u-neq-0}
\end{equation}
which allows to consider the Dirichlet quotients for every $t\geq 0$.
In the second section of the proof we assume that there exist a
constant $c_{1}>0$ and a sequence $t_{n}\to +\infty$ such that
\begin{equation}
	|A^{1/2}u(t_{n})|^{2}\leq c_{1}|u(t_{n})|^{2+p}
	\quad\quad
	\forall n\in\n.
	\label{hp:liminf}
\end{equation}

Under this assumption, we prove that a similar estimate holds true 
for all times, namely
\begin{equation}
	|A^{1/2}u(t)|^{2}\leq c_{2}|u(t)|^{2+p}
	\quad\quad
	\forall t\geq 0
	\label{hp:liminf-glob}
\end{equation}
for a suitable constant $c_{2}\geq c_{1}$.  This is not yet
(\ref{th:range}), but in any case it shows that $|A^{1/2}u(t)|$ decays
faster than $|u(t)|$.  As already pointed out, this means that the
solution $u(t)$ moves closer and closer to the kernel of $A$, and
suggests that the terms with $Au(t)$ and $A^{1/2}u(t)$ in
equation~(\ref{pbm:linear}) and estimate~(\ref{hp:diff-ineq}) can be
neglected.  With this ansatz, we obtain that $|u'(t)|\leq
K_{0}|u(t)|^{1+p}$, and it is easy to show that all nonzero solutions
of this differential inequality are slow in the sense of
(\ref{th:slow}).  Finally, we improve (\ref{hp:liminf-glob}) in order
to obtain (\ref{th:range}).

In the third and last section of the proof we are left with the case where
(\ref{hp:liminf}) is false for every constant $c_{1}$ and every
sequence $t_{n}\to+\infty$.  This easily implies that there exists
$T_{0}\geq 0$ such that
\begin{equation}
	|u(t)|^{2+p}\leq|A^{1/2}u(t)|^{2}
	\quad\quad
	\forall t\geq T_{0}.
	\label{est:u<Au}
\end{equation}

This means that now $u(t)$ is faraway from the kernel of $A$.  Thus we
are not allowed to ignore the operator $A$, but we can neglect the
right-hand side of (\ref{pbm:linear}) because the exponents in
(\ref{hp:diff-ineq}) are larger than one.  Therefore, a good
approximation of (\ref{pbm:linear}) is now the linear homogeneous
equation $u'(t)+Au(t)=0$, whose solutions decay exponentially with
possible rates corresponding to eigenvalues of $A$.  The formal proof
requires several steps.  First of all, we provide exponential
estimates from below and from above with non-optimal rates.  Then we
identify the exact rate, and finally we prove that (\ref{th:fast+})
holds true.

We point out that the exponent $2+p$ is non-optimal both in
(\ref{hp:liminf}) and in the opposite estimate (\ref{est:u<Au}).
Indeed, a posteriori it turns out that (up to multiplicative
constants) $|A^{1/2}u|\leq|u|^{1+p}$ in the case of slow
solutions, and $|A^{1/2}u|\sim|u|$ in the case of fast solutions, so
that the right exponents would be $2+2p$ and $2$, respectively.
Nevertheless, the intermediate exponent $2+p$ acts as a threshold
separating the two different regimes, and leaving enough room on both
sides to perform our estimates.

\subsubsection*{Non-trivial solutions never vanish.}

\paragraph{\textmd{\textit{Forward uniqueness}}}

We prove that $u(0)=0$ implies that $u(t)=0$ for every $t\geq 0$.

To this end, we set
$z(t):=|u(t)|^{2}+|A^{1/2}u(t)|^{2}$.  A simple computation shows that
\begin{eqnarray*}
	z'(t) & = & -2|A^{1/2}u(t)|^{2}-2|Au(t)|^{2}
	+2\langle u(t),g(t)\rangle+2\langle Au(t),g(t)\rangle\\
	 & \leq & -2|A^{1/2}u(t)|^{2}-2|Au(t)|^{2}+|u(t)|^{2}+|g(t)|^{2}
	 +|Au(t)|^{2}+|g(t)|^{2}\\
	 & \leq & |u(t)|^{2}+2|g(t)|^{2}.
\end{eqnarray*}

From (\ref{hp:diff-ineq}) we obtain that
\begin{equation}
	|g(t)|^{2}\leq 2K_{0}^{2}
	\left(|u(t)|^{2+2p}+|A^{1/2}u(t)|^{2+2q}\right),
	\label{est:gt2}
\end{equation}
hence
$$z'(t)\leq z(t)+4K_{0}^{2}[z(t)]^{1+p}+
4K_{0}^{2}[z(t)]^{1+q}
\quad\quad
\forall t\geq 0.$$

All powers of $z(t)$ in the right-hand side of this scalar 
differential inequality have exponents greater
than or equal to one.  It follows that the right-hand side, as a
function of $z(t)$, is Lipschitz continuous.  This is enough to
guarantee that necessarily $z(t)=0$, hence $u(t)=0$, for every $t\geq
0$.

\paragraph{\textmd{\textit{Backward uniqueness}}}

We prove that $u(0)\neq 0$ implies that $u(t)\neq 0$ for every $t\geq
0$.

As in the classical references~\cite{bt,ghidaglia}, here we exploit
the standard Dirichlet quotient $Q(t)$ defined in (\ref{defn:Q}).  To
this end, we set
$$S:=\sup\left\{t\geq 0:u(\tau)\neq 0\quad\forall\tau\in[0,t]\right\}.$$

Since $u(0)\neq 0$, a simple continuity argument shows that $S>0$. 
Backward uniqueness is equivalent to saying that $S=+\infty$. So let us 
assume by contradiction that $S<+\infty$. By the maximality of $S$, 
this means that $u(S)=0$. Now we show that this is not possible.

In the interval $[0,S)$ we have that $u(t)\neq 0$, hence $Q(t)$ is
defined.  Let us estimate its time-derivative as in
(\ref{th:deriv-Q}).  If we neglect the first term in the right-hand
side, and we estimate the second one by means of (\ref{est:gt2}), we
obtain that
\begin{eqnarray*}
	Q'(t) & \leq & \frac{1}{|u(t)|^{2}}\cdot
	2K_{0}^{2}\left(|u(t)|^{2+2p}+|A^{1/2}u(t)|^{2+2q}\right)
	\\
	\noalign{\vspace{1ex}}
	 & \leq &
	2K_{0}^{2}|u(t)|^{2p}+2K_{0}^{2}|A^{1/2}u(t)|^{2q}\cdot Q(t).
\end{eqnarray*}

Due to the regularity of the solution, and in particular
to~(\ref{th:reg-cont}), both $|u(t)|$ and $|A^{1/2}u(t)|$ are bounded
on bounded time-intervals, and in particular in $[0,S)$.  This is
enough to conclude that also $Q(t)$ is bounded in $[0,S)$.

Now let us consider the function $y(t):=|u(t)|^{2}$. A simple 
computation shows that
\begin{equation}
	y'(t)=-2|A^{1/2}u(t)|^{2}+2\langle u(t),g(t)\rangle\geq
	-2|A^{1/2}u(t)|^{2}-2|u(t)|\cdot|g(t)|.
	\label{eqn:y'}
\end{equation}

By definition of Dirichlet quotient, assumption (\ref{hp:diff-ineq}) 
can be rewritten as
$$|g(t)|\leq  K_{0}|u(t)|\left(|u(t)|^{p}+
[Q(t)]^{1/2}|A^{1/2}u(t)|^{q}\right),$$
so that
\begin{eqnarray}
	y'(t) & \geq & -2\frac{|A^{1/2}u(t)|^{2}}{|u(t)|^{2}}|u(t)|^{2}-
	2K_{0}|u(t)|^{2}\left(|u(t)|^{p}+
	[Q(t)]^{1/2}|A^{1/2}u(t)|^{q}\strut\right)
	\nonumber  \\
	\noalign{\vspace{1ex}}
	 & = & -2\left(Q(t)+K_{0}|u(t)|^{p}+
	 K_{0}[Q(t)]^{1/2}|A^{1/2}u(t)|^{q}\right)y(t).
	 \label{est:y'}
\end{eqnarray}

Since $|u(t)|$, $|A^{1/2}u(t)|$, and $Q(t)$ are bounded in $[0,S)$, we
deduce that there exists $c_{3}$ such that $y'(t)\geq -c_{3}y(t)$ for
almost every $t\in[0,S)$, hence 
$$y(t)\geq y(0)e^{-c_{3}t} 
\quad\quad
\forall t\in[0,S).$$

Since $y(0)> 0$, letting $t\to S^{-}$ we conclude that $y(S)> 0$,
hence $u(S)\neq 0$.  This contradicts the maximality of $S$, and
completes the proof that $u(t)$ cannot vanish in a finite time.

\subsubsection*{Slow solutions}

In this second part of the proof we consider the case where $u(t)$ in
not the null solution and (\ref{hp:liminf}) holds true for some
$c_{1}> 0$ and some sequence $t_{n}\to +\infty$.

\paragraph{\textmd{\textit{Main estimate}}}

We prove that there exists a constant $c_{2}$ such that
(\ref{hp:liminf-glob}) holds true.  

This estimate is trivial if $A$ is the null operator.  Otherwise, let
$\nu>0$ denote the smallest positive eigenvalue of $A$, which exists
because we assumed that eigenvalues are an increasing sequence.  With 
this choice, the operator $A$ satisfies assumption~(\ref{hp:A-nu}).  Let
us consider the modified Dirichlet quotient~(\ref{defn:Qd}) with
$d:=p$, which is defined for every $t\geq 0$ by virtue of
(\ref{hp:u-neq-0}).  From (\ref{th:deriv-Qd}) it follows that
$$Q_{p}'(t)\leq -\nu Q_{p}(t)+ 2(2+p)|u(t)|^{p}\cdot[Q_{p}(t)]^{2}+
(3+p)\frac{|g(t)|^{2}}{|u(t)|^{2+p}}.$$
	
Therefore, if we write (\ref{est:gt2}) in the form
$$|g(t)|^{2}\leq 2K_{0}^{2}|u(t)|^{2+p}\left(
|u(t)|^{p}+[Q_{p}(t)]^{1+q}|u(t)|^{(2+p)q}\right),$$
we obtain that
\begin{eqnarray}
	Q_{p}'(t) & \leq & -\nu Q_{p}(t)+ 
	2(2+p)|u(t)|^{p}\cdot[Q_{p}(t)]^{2}
	\nonumber
	\\
	\noalign{\vspace{1ex}}
	 &  & +2(3+p)K_{0}^{2}|u(t)|^{p}+
	 2(3+p)K_{0}^{2}\cdot[Q_{p}(t)]^{1+q}|u(t)|^{(2+p)q}.
	 \label{est:Qp'-good}
\end{eqnarray}

Let us consider now the constant $c_{1}$ and the sequence $t_{n}\to
+\infty$ of (\ref{hp:liminf}).  From assumption~(\ref{hp:u-limit}) we
have in particular that $|u(t)|\to 0$ as $t\to +\infty$.  Therefore,
there exists $n_{0}\in\n$ such that 
\begin{eqnarray}
	2\nu c_{1} & \geq & 2(2+p)|u(t)|^{p}\cdot(2c_{1})^{2}
	+2(3+p)K_{0}^{2}|u(t)|^{p}
	\nonumber  \\
	\noalign{\vspace{1ex}}
	 &  & +2(3+p)K_{0}^{2}\cdot(2c_{1})^{1+q}|u(t)|^{(2+p)q}
	\label{defn:tn0}
\end{eqnarray}
for every $t\geq t_{n_{0}}$. Now we claim that
\begin{equation}
	Q_{p}(t)\leq 2c_{1}
	\quad\quad
	\forall t\geq t_{n_{0}}.
	\label{est:Qp}
\end{equation}

Since $Q_{p}(t)$ is continuous, and consequently bounded, in the
compact interval $[0,t_{n_{0}}]$, this is enough to establish
(\ref{hp:liminf-glob}).

In order to prove (\ref{est:Qp}), we set
$$S:=\sup\left\{t\geq t_{n_{0}}:Q_{p}(\tau)\leq 2c_{1}
\quad\forall\tau\in[t_{n_{0}},t]\right\},$$
so that now (\ref{est:Qp}) is equivalent to $S=+\infty$.  To begin
with, we observe that $Q_{p}(t_{n_{0}})\leq c_{1}<2c_{1}$, hence a
simple continuity argument yields that $S>t_{n_{0}}$.  Let us assume
by contradiction that $S<+\infty$.  By the maximality of $S$, this
means that $Q_{p}(S)=2c_{1}$.  

Now we show that this is impossible.  We already know that
$Q_{p}(t)\leq 2c_{1}$ for every $t\in[t_{n_{0}},S)$.  Plugging
this estimate into (\ref{est:Qp'-good}), and exploiting
(\ref{defn:tn0}), we obtain that
$$Q_{p}'(t)\leq -\nu Q_{p}(t)+2\nu c_{1}
\quad\quad
\forall t\in[t_{n_{0}},S).$$

Integrating this differential inequality, and recalling once again
that $Q_{p}(t_{n_{0}})\leq c_{1}$, we obtain that
$$Q_{p}(t)\leq 2c_{1}-c_{1}\exp\left(-\nu (t-t_{n_{0}})\right)
\quad\quad
\forall t\in[t_{n_{0}},S).$$

Letting $t\to S^{-}$, we conclude that $Q_{p}(S)<2c_{1}$, which
contradicts the maximality of $S$.  This completes the proof of
(\ref{hp:liminf-glob}).

\paragraph{\textmd{\textit{Faster decay of the range component}}}

Let us prove (\ref{th:range}).  Once again this is trivial if $A$ is
the null operator.  Otherwise, assumption~(\ref{hp:A-nu}) is satisfied
with $\nu>0$ equal to the smallest positive eigenvalue of $A$.  Let us
consider the modified Dirichlet quotient~(\ref{defn:Qd}) with $d:=2p$,
which is defined for every $t\geq 0$ because of (\ref{hp:u-neq-0}).
From (\ref{th:deriv-Qd}) it follows that
\begin{equation}
	Q_{2p}'(t)\leq -\nu Q_{2p}(t)+
	4(1+p)|u(t)|^{2p}\cdot[Q_{2p}(t)]^{2}+
	(3+2p)\frac{|g(t)|^{2}}{|u(t)|^{2+2p}}.	
	\label{est:Q2p'}
\end{equation}

From (\ref{hp:liminf-glob}) we have that
\begin{equation}
	|u(t)|^{2p}\cdot[Q_{2p}(t)]^{2}=
	|u(t)|^{p}\cdot Q_{2p}(t)\cdot Q_{p}(t)\leq
	c_{2}|u(t)|^{p}\cdot Q_{2p}(t).
	\label{est:Q2p'-1}
\end{equation}

Moreover, now we can rewrite (\ref{est:gt2}) in the form 
\begin{equation}
	|g(t)|^{2}\leq 2K_{0}^{2}|u(t)|^{2+2p}\left(
	1+[Q_{2p}(t)]\cdot|A^{1/2}u(t)|^{2q}\right).
	\label{est:Q2p'-2}
\end{equation}

Plugging (\ref{est:Q2p'-1}) and (\ref{est:Q2p'-2}) into 
(\ref{est:Q2p'}), we obtain that
$$Q_{2p}'(t)\leq -Q_{2p}(t)\cdot\left\{
\nu-c_{4}|u(t)|^{p}-c_{5}|A^{1/2}u(t)|^{2q}
\right\}+c_{6}.$$

Due to assumption~(\ref{hp:u-limit}), there exists $T_{1}\geq 0$ such
that
$$\nu-c_{4}|u(t)|^{p}-c_{5}|A^{1/2}u(t)|^{2q}
\geq\frac{\nu}{2}
\quad\quad
\forall t\geq T_{1},$$
hence
$$Q_{2p}'(t)\leq -\frac{\nu}{2}Q_{2p}(t)+c_{6}$$
for almost every $t\geq T_{1}$.  Integrating this differential
inequality we conclude that $Q_{2p}(t)$ is uniformly bounded for every
$t\geq T_{1}$.  Since $Q_{2p}(t)$ is continuous, and consequently
bounded, in the compact interval $[0,T_{1}]$, this is enough to
establish (\ref{th:range}).

\paragraph{\textmd{\textit{Slow decay of the solution}}}

Let us set as usual $y(t):=|u(t)|^{2}$.  Since $Q_{2p}(t)$ and
$|A^{1/2}u(t)|$ are uniformly bounded, from (\ref{est:Q2p'-2}) we have
now that
$$|g(t)|\leq c_{7}|u(t)|^{1+p}
\quad\quad
\forall t\geq 0.$$

Plugging this estimate and (\ref{hp:liminf-glob}) into (\ref{eqn:y'})
we obtain that
$$y'(t)\geq -2|A^{1/2}u(t)|^{2}-2|g(t)|\cdot|u(t)|\geq 
-c_{8}|u(t)|^{2+p}=-c_{8}[y(t)]^{1+p/2}.$$

Integrating this differential inequality we deduce (\ref{th:slow}).

\subsubsection*{Spectral fast solutions}

In this last section of the proof it remains to consider the case
where (\ref{hp:u-neq-0}) and (\ref{est:u<Au}) hold true.  This implies
in particular that $|A^{1/2}u(t)|\neq 0$ for every $t\geq T_{0}$, and
therefore $A$ is not the null operator.  

\paragraph{\textmd{\textit{Non-optimal exponential decay from above}}}

Let $\nu$ be the smallest positive eigenvalue of $A$.
We prove that there exists a constant
$c_{9}$ such that
\begin{equation}
	|A^{1/2}u(t)|\leq c_{9}\exp\left(-\frac{\nu}{4}t\right)
	\quad\quad
	\forall t\geq 0.
	\label{th:Au-exp-above}
\end{equation}

As a consequence of (\ref{est:u<Au}), this implies also that
\begin{equation}
	|u(t)|\leq c_{10}\exp\left(-\frac{\nu}{2(2+p)}t\right)
	\quad\quad
	\forall t\geq 0.
	\label{th:u-exp-above}
\end{equation}

To this end, we consider the function $E(t):=|A^{1/2}u(t)|^{2}$, and 
we estimate its time-derivative as usual
\begin{equation}
	E'(t)=-2|Au(t)|^{2}+2\langle Au(t),g(t)\rangle\leq
	-|Au(t)|^{2}+|g(t)|^{2}.
	\label{est:E'}
\end{equation}

For the first term we have that 
$$-|Au(t)|^{2}\leq -\nu|A^{1/2}u(t)|^{2}=-\nu E(t).$$

For the second term, from (\ref{est:gt2}) and (\ref{est:u<Au}) it 
follows that
\begin{equation}
	|g(t)|^{2}\leq 
	2K_{0}^{2}|A^{1/2}u(t)|^{2}\left(|u(t)|^{p}+|A^{1/2}u(t)|^{2q}\right)
	\quad\quad
	\forall t\geq T_{0}.
	\label{est:gt-fast}
\end{equation}

Plugging these estimates into (\ref{est:E'}) we obtain that
$$E'(t)\leq -\left(
\nu-2K_{0}^{2}|u(t)|^{p}-2K_{0}^{2}|A^{1/2}u(t)|^{2q}
\right)E(t).$$

Exploiting again assumption (\ref{hp:u-limit}), we deduce that
there exists $T_{1}\geq T_{0}$ such that
$$\nu-2K_{0}^{2}|u(t)|^{p}-2K_{0}^{2}|A^{1/2}u(t)|^{2q}\geq
\frac{\nu}{2}
\quad\quad
\forall t\geq T_{1}.$$

It follows that $E'(t)\leq -(\nu/2) E(t)$ for every $t\geq T_{1}$,
hence
$$E(t)\leq E(T_{1})\exp\left(-\frac{\nu}{2}(t-T_{1})\right)
\quad\quad
\forall t\geq T_{1},$$
which easily implies (\ref{th:Au-exp-above}).

\paragraph{\textmd{\textit{Boundedness of the Dirichlet quotient}}}

Let us consider once again the Dirichlet quotient $Q(t)$ defined by
(\ref{defn:Q}), which now is defined for every $t\geq 0$ because of 
(\ref{hp:u-neq-0}). We prove that there exists a constant $c_{11}$ 
such that
\begin{equation}
	Q(t)\leq c_{11}
	\quad\quad
	\forall t\geq 0.
	\label{th:Q-bounded}
\end{equation}

To this end, we estimate $Q'(t)$ starting from (\ref{th:deriv-Q}).
Exploiting (\ref{est:gt-fast}) we obtain that
$$Q'(t)\leq\frac{|g(t)|^{2}}{|u(t)|^{2}}\leq
2K_{0}^{2}Q(t)
\left(|u(t)|^{p}+|A^{1/2}u(t)|^{2q}\right)$$
for almost every $t\geq T_{0}$.  Integrating this differential
inequality we find that
$$Q(t)\leq Q(T_{0})\exp\left(2K_{0}^{2}\int_{T_{0}}^{t}
\left(|u(s)|^{p}+|A^{1/2}u(s)|^{2q}\right)\,ds\right)
\quad\quad
\forall t\geq T_{0}.$$

Due to (\ref{th:u-exp-above}) and (\ref{th:Au-exp-above}), the
integral in the right-hand side is bounded independently of $t$.
Since $Q(t)$ is continuous, and consequently bounded, in the compact
interval $[0,T_{0}]$, this is enough to establish
(\ref{th:Q-bounded}).

\paragraph{\textmd{\textit{Non-optimal exponential decay from below}}}

We prove that there exist positive constants
$c_{12}$ and $c_{13}$ such that
\begin{equation}
	|u(t)|\geq c_{12}e^{-c_{13}t}
	\quad\quad
	\forall t\geq 0.
	\label{th:u-exp-below}
\end{equation}

To this end, we consider once again the function $y(t):=|u(t)|^{2}$,
and we estimate $y'(t)$ starting from (\ref{est:y'}).  Now we have
that $|u(t)|$ and $|A^{1/2}u(t)|$ are bounded independently of $t$
because of assumption~(\ref{hp:u-limit}), and $Q(t)$ is bounded
independently of $t$ because of (\ref{th:Q-bounded}).  Therefore,
there exists a constant $c_{14}$ such that
$y'(t)\geq -c_{14}y(t)$ for almost every $t\geq 0$.  Since $y(0)> 0$,
integrating this differential inequality we obtain
(\ref{th:u-exp-below}).
 
\paragraph{\textmd{\textit{Exact exponential decay rate}}}

Let us set
\begin{equation}
	\lambda:=\sup\left\{\gamma>0:
	\lim_{t\to+\infty}|u(t)|_{D(A^{1/2})}e^{\gamma t}=0\right\}.
	\label{defn:lambda}
\end{equation}

From (\ref{th:Au-exp-above}) and (\ref{th:u-exp-above}) it follows
that $\lambda$ is the supremum of a nonempty set.  Moreover, from
(\ref{th:u-exp-below}) it follows that $|u(t)|_{D(A^{1/2})}\geq
c_{12}e^{-c_{13}t}$, which implies that $\lambda$ is finite.
Therefore, $\lambda$ is a positive real number.

We claim that $\lambda$ is an eigenvalue of $A$.  To this end, we
write $H$ as an orthogonal direct sum
\begin{equation}
   H:=H_{\lambda,-}\oplus H_{\lambda}\oplus H_{\lambda,+},
   \label{defn:H-lambda}
\end{equation}
where 
\begin{itemize}
	\item  $H_{\lambda}$ is the eigenspace relative to $\lambda$ if 
	$\lambda$ is an eigenvalue of $A$, or $H_{\lambda}=\{0\}$ 
	otherwise,

	\item $H_{\lambda,-}$ is the closure of the space generated by all
	eigenvectors relative to eigenvalues of $A$ less than $\lambda$
	(if any),

	\item $H_{\lambda,+}$ is the closure of the space generated by all
	eigenvectors relative to eigenvalues of $A$ greater than $\lambda$
	(if any).
\end{itemize}

These three subspaces of $H$ are $A$-invariant, and some of them might be the
trivial subspace $\{0\}$ depending on the value of $\lambda$.  Let
$u_{\lambda,-}(t)$, $u_{\lambda}(t)$ and $u_{\lambda,+}(t)$ denote the
components of $u(t)$ with respect to the decomposition
(\ref{defn:H-lambda}), and let $g_{\lambda,-}(t)$, $g_{\lambda}(t)$
and $g_{\lambda,+}(t)$ be the corresponding components of $g(t)$.  Let
$\beta$ be the smallest eigenvalue of $A$ larger than $\lambda$ (if
any), or $\beta=+\infty$ otherwise, and let 
$\delta:=\min\{(1+p)\lambda,(1+q)\lambda\}$. 

First of all, we observe that our definition of $\lambda$, combined 
with assumption (\ref{hp:diff-ineq}), implies that
\begin{equation}
	\lim_{t\to +\infty}|g(t)|e^{\gamma t}=0
	\quad\quad
	\forall\gamma<\delta.
	\label{est:gt-exp}
\end{equation}

We claim that 
\begin{equation}
	\lim_{t\to +\infty}
	|u_{\lambda,+}(t)|_{D(A^{1/2})}e^{\gamma t}=0
	\quad\quad
	\forall\gamma<\min\{\beta,\delta\},
	\label{est:u+}
\end{equation}
\begin{equation}
	\lim_{t\to +\infty}
	|u_{\lambda,-}(t)|_{D(A^{1/2})}e^{\gamma t}=0
	\quad\quad
	\forall\gamma<\delta.
	\label{est:u-}
\end{equation}

This is enough to conclude that $\lambda$ is an eigenvalue of $A$,
because otherwise $u_{\lambda}(t)\equiv 0$, so that (\ref{est:u+}) and
(\ref{est:u-}) would imply that
\begin{equation}
	\lim_{t\to +\infty}|u(t)|_{D(A^{1/2})}e^{\gamma t}=0
	\quad\quad
	\forall\gamma<\min\{\beta,\delta\},
	\label{est:u-glob}
\end{equation}
and this would contradict the maximality of $\lambda$ because
$\min\{\beta,\delta\}>\lambda$.

Let us prove (\ref{est:u+}).  If all eigenvalues of $A$
are less than or equal to $\lambda$, then $H_{\lambda,+}=\{0\}$, so
that (\ref{est:u+}) is trivial.  Otherwise, we can apply
Lemma~\ref{lemma:ODE-super} with 
$$X:=H_{\lambda,+},
\hspace{3em}
w(t):=u_{\lambda,+}(t),
\hspace{3em}
\psi(t):=g_{\lambda,+}(t).$$

Indeed assumptions (\ref{hp:ODE-spA}) and (\ref{hp:ODE-sppsi}) follow
from our definition of $\beta$ and from estimate (\ref{est:gt-exp}).
At this point, conclusion (\ref{th:ODE-sp}) of
Lemma~\ref{lemma:ODE-super} is exactly (\ref{est:u+}).

Let us prove (\ref{est:u-}).  If all eigenvalues of $A$ are
greater than or equal to $\lambda$, then $H_{\lambda,-}=\{0\}$, so
that (\ref{est:u-}) is trivial.  Otherwise, let $\alpha$ be the
largest eigenvalue of $A$ less than $\lambda$.  In this case we can
apply Lemma~\ref{lemma:ODE-sub} with 
$$X:=H_{\lambda,-},
\hspace{3em}
w(t):=u_{\lambda,-}(t),
\hspace{3em}
\psi(t):=g_{\lambda,-}(t).$$

Indeed assumptions (\ref{hp:ODE-sbA}) and (\ref{hp:ODE-sbpsi}) follow
from our definition of $\alpha$ and from estimate (\ref{est:gt-exp}).
From Lemma~\ref{lemma:ODE-sub} we deduce that $e^{tA}u_{\lambda,-}(t)$
has a limit $x_{0}\in H_{\lambda,-}$ as $t\to +\infty$.  Recalling
that $\alpha<\lambda$, from (\ref{defn:lambda}) we have that
$$|x_{0}|=\lim_{t\to+\infty}\left|e^{tA}u_{\lambda,-}(t)\right|\leq
\lim_{t\to+\infty}e^{\alpha t}|u_{\lambda,-}(t)|\leq
\lim_{t\to+\infty}|u(t)|e^{\alpha t}=0,$$
so that $x_{0}=0$.  At this point conclusion (\ref{th:remainder}) of
Lemma~\ref{lemma:ODE-sub} holds true with $x_{0}=0$, and this is 
exactly to (\ref{est:u-}).

\paragraph{\textmd{\textit{Exact limit}}}

Now we know that $\lambda$ is an eigenvalue of $A$, and that
$u_{\lambda,-}(t)$ and $u_{\lambda,+}(t)$ decay faster than
$e^{-\gamma t}$ for every $\gamma$ satisfying (\ref{defn:delta}).  It
remains to consider the component $u_{\lambda}(t)$.  To this end, we
apply again Lemma~\ref{lemma:ODE-sub}, this time with
$$X:=H_{\lambda},
\hspace{3em}
w(t):=u_{\lambda}(t),
\hspace{3em}
\psi(t):=g_{\lambda}(t).$$  

Now assumption (\ref{hp:ODE-sbA}) is trivially satisfied with
$\alpha:=\lambda$ because the operator $A$ is $\lambda$ times the
identity in $H_{\lambda}$, while assumption (\ref{hp:ODE-sbpsi})
follows again from estimate (\ref{est:gt-exp}).

Thus from Lemma~\ref{lemma:ODE-sub} we deduce that 
$e^{\lambda t}u_{\lambda}(t)$ tends to some $v_{0}\in H_{\lambda}$ as 
$t\to +\infty$, and 
\begin{equation}
	\lim_{t\to +\infty}
	\left|u_{\lambda}(t)-v_{0}e^{-\lambda t}\right|e^{\gamma t}=0
	\quad\quad
	\forall\gamma<\delta.
	\label{th:limit-lambda}
\end{equation}

We claim that $v_{0}\neq 0$.  Indeed otherwise
(\ref{th:limit-lambda}), (\ref{est:u+}) and (\ref{est:u-}) would imply
(\ref{est:u-glob}), and this would contradict the maximality of
$\lambda$ because $\min\{\beta,\delta\}>\lambda$.

At this point (\ref{est:u+}), (\ref{est:u-}) and
(\ref{th:limit-lambda}) imply (\ref{th:fast+}) for every
$\gamma<\min\{\beta,\delta\}$, hence for every $\gamma$ satisfying
(\ref{defn:delta}).\qed

\subsection{Proof of Theorem~\ref{thm:main-slow}}

Let us set 
$$K_{1}=\frac{4K_{0}^{2}(3+2p)}{\nu},$$
and let us choose $\sigma_{0}>0$ small enough so that the following 
two conditions are satisfied:
\begin{equation}
	\sigma_{0}^{2}+K_{1}\sigma_{0}^{2+2p}<R^{2},
	\label{hp:sigma-0-1}
\end{equation}
\begin{equation}
	4(1+p)\sigma_{0}^{2p}K_{1}^{2}+
	2K_{0}^{2}(3+2p)K_{1}^{1+q}\sigma_{0}^{(2+2p)q}\leq
	2K_{0}^{2}(3+2p).
	\label{hp:sigma-0-2}
\end{equation}

Let $\mathcal{S}$ be the set of all $u_{0}\in D(A^{1/2})$ such that
\begin{equation}
	u_{0}\neq 0,
	\quad\quad\quad
	|u_{0}|<\sigma_{0},
	\quad\quad\quad
	|A^{1/2}u_{0}|^{2}<K_{1}|u_{0}|^{2+2p}.
	\label{hp:u0}
\end{equation}

It is clear that these assumptions define an open set in $D(A^{1/2})$,
which is nonempty because it contains at least all $u_{0}\in\ker(A)$
with $u_{0}\neq 0$ and $|u_{0}|<\sigma_{0}$.

Let $u_{0}\in\mathcal{S}$, and let $u(t)$ be the unique local
solution to problem (\ref{pbm:nl})--(\ref{pbm:data}) provided by
Theorem~\ref{thmbibl:nl}, defined in a maximal interval $[0,T)$.  We
claim that $T=+\infty$ and this solution is slow in the sense of
(\ref{th:slow}).

\subparagraph{\textmd{\textit{Basic estimate}}}

We prove that
\begin{equation}
	|u(t)|\leq|u_{0}|<\sigma_{0}
	\quad\quad
	\forall t\in[0,T).
	\label{th:ut<u0}
\end{equation}

Indeed let us set $y(t):=|u(t)|^{2}$. A simple computation shows that
\begin{equation}
	y'(t)=-2|A^{1/2}u(t)|^{2}+2\langle u(t),f(u(t))\rangle.
	\label{eqn:y'-bis}
\end{equation}

Thanks to~(\ref{hp:f-sign}) we have that $y'(t)\leq 0$, which proves
(\ref{th:ut<u0}).

\subparagraph{\textmd{\textit{Boundedness of the generalized Dirichlet 
quotient}}}

We prove that 
\begin{equation}
	u(t)\neq 0
	\quad\mbox{and}\quad
	|A^{1/2}u(t)|^{2}<K_{1}|u(t)|^{2+2p}
	\quad\quad
	\forall t\in[0,T).
	\label{est:up-to-T}
\end{equation}

To this end, we consider once again the modified Dirichlet quotient 
(\ref{defn:Qd}) with $d:=2p$, defined as long as $u(t)\neq 0$, and we set
$$S:=\sup\left\{t\in[0,T):
|u(\tau)|\cdot(Q_{2p}(\tau)-K_{1})<0
\quad\forall\tau\in[0,t]\right\},$$
so that now (\ref{est:up-to-T}) is equivalent to $S=T$.  To begin
with, we set $t=0$ and we obtain that $|u(0)|\cdot(Q_{2p}(0)-K_{1})<0$
because of the first and third condition in~(\ref{hp:u0}).  Therefore,
a simple continuity argument gives that $S>0$.  Let us assume by
contradiction that $S<T$.  Then, by the maximality of $S$, this means
that at least one of the following equalities is satisfied
\begin{equation}
	\quad
	u(S)=0,
	\hspace{3em}
	Q_{2p}(S)=K_{1}. 
	\label{eqn:alternative}
\end{equation}

Now we exclude both possibilities. First of all, our 
definition of $S$ implies that
\begin{equation}
	u(t)\neq 0
	\quad\mbox{and}\quad
	Q_{2p}(t)<K_{1}
	\quad\quad
	\forall t\in[0,S).
	\label{est:up-to-S}
\end{equation}

In particular, keeping into account assumption (\ref{hp:f-ho}) and 
(\ref{th:ut<u0}), we obtain that
\begin{eqnarray}
	|f(u(t))| & \leq & K_{0}|u(t)|^{1+p}\left(
	1+[Q_{2p}(t)]^{(1+q)/2}|u(t)|^{(1+p)q}\right)
	\nonumber  \\
	\noalign{\vspace{1ex}}
	 & \leq & K_{0}|u(t)|^{1+p}\left(
	1+K_{1}^{(1+q)/2}\sigma_{0}^{(1+p)q}\right)
	\label{est:f-new}
\end{eqnarray}
for every $t\in[0,S)$.  Let us consider again the function
$y(t):=|u(t)|^{2}$.  We compute the time-derivative as in
(\ref{eqn:y'-bis}), and then we estimate the right-hand side
exploiting (\ref{est:up-to-S}) and (\ref{est:f-new}).  We obtain that
\begin{eqnarray}
	y'(t) & \geq & -2|A^{1/2}u(t)|^{2}
	-2|u(t)|\cdot|f(u(t))|  
	\nonumber \\
	\noalign{\vspace{1ex}}
	 & \geq & -2Q_{2p}(t)\cdot|u(t)|^{p}\cdot|u(t)|^{2+p}-
	 2K_{0}|u(t)|^{2+p}\left(
	 1+K_{1}^{(1+q)/2}\sigma_{0}^{(1+p)q}\right)   
	 \nonumber  \\
	\noalign{\vspace{1ex}}
	 & \geq & -2\left(
	 K_{1}\sigma_{0}^{p}+K_{0}+K_{0}K_{1}^{(1+q)/2}\sigma_{0}^{(1+p)q}
	 \right)|y(t)|^{1+p/2}
	 \label{est:y'-slow}
\end{eqnarray}
for almost every $t\in[0,S)$.  Integrating this differential
inequality we conclude that $u(t)$ cannot vanish in a finite time,
which rules out the first possibility in~(\ref{eqn:alternative}).

In order to exclude the second one, we compute the time-derivative of
$Q_{2p}(t)$. From (\ref{th:deriv-Qd}) with $d=2p$ we obtain that
$$Q_{2p}'(t) \leq -\nu Q_{2p}(t)+ 4(1+p)|u(t)|^{2p}[Q_{2p}(t)]^{2}+
(3+2p)\frac{|f(u(t))|^{2}}{|u(t)|^{2+2p}}.$$

From (\ref{th:ut<u0}), (\ref{est:up-to-S}) and (\ref{est:f-new}) we have that 
$$|u(t)|^{2p}[Q_{2p}(t)]^{2}\leq \sigma_{0}^{2p}K_{1}^{2},
\hspace{3em}
\frac{|f(u(t))|^{2}}{|u(t)|^{2+2p}}\leq 2K_{0}^{2}\left(
1+K_{1}^{1+q}\sigma_{0}^{(2+2p)q}\right),$$
so that
$$Q_{2p}'(t)\leq -\nu Q_{2p}(t)+4(1+p)\sigma_{0}^{2p}K_{1}^{2}+
2K_{0}^{2}(3+2p)+2K_{0}^{2}(3+2p)K_{1}^{1+q}\sigma_{0}^{(2+2p)q}.$$

Keeping the smallness condition (\ref{hp:sigma-0-2}) into account, we 
finally deduce that
$$Q_{2p}'(t)\leq -\nu Q_{2p}(t)+4K_{0}^{2}(3+2p)=
-\nu(Q_{2p}(t)-K_{1})$$
for almost every $t\in[0,S)$.  Integrating this differential
inequality we obtain that
$$Q_{2p}(t)\leq K_{1}+(Q_{2p}(0)-K_{1})e^{-\nu t}
\quad\quad
\forall t\in[0,S).$$

Letting $t\to S^{-}$, and recalling that $Q_{2p}(0)<K_{1}$, we
conclude that $Q_{2p}(S)<K_{1}$, which rules out the second
possibility in~(\ref{eqn:alternative}).

\subparagraph{\textmd{\textit{Global existence and slow decay}}}

We show that $T=+\infty$.  Let us assume indeed that $T<+\infty$.
Letting $t\to T^{-}$ in (\ref{th:ut<u0}) and (\ref{est:up-to-T}), and
taking into account the smallness assumption (\ref{hp:sigma-0-1}), we
obtain that
$$\lim_{t\to T^{-}}\left(|u(t)|^{2}+|A^{1/2}u(t)|^{2}\right)\leq
\sigma_{0}^{2}+K_{1}\sigma_{0}^{2+2p}<R^{2},$$ 
which contradicts (\ref{th:alternative}).

Since we have proved that $S=T=+\infty$, the differential inequality in 
(\ref{est:y'-slow}) now holds true for every $t\geq 0$. A simple 
integration of this differential inequality proves that $u(t)$
satisfies (\ref{th:slow}) for a suitable positive constant $M_{1}$, 
depending only on $K_{0}$, $K_{1}$, $\sigma_{0}$, $p$, $q$.\qed

\subsection{Proof of Theorem~\ref{thm:main-exponential}}

Let us sketch the strategy of the proof, based on a fixed point
argument, before entering into details.  We begin by observing that 
our definition of $H_{-}$ implies that $H_{-}\subseteq D(A)$ 
and 
\begin{equation}
	|Au|\leq\lambda|u|
	\quad\quad
	\forall u\in H_{-},
   \label{hp:A-sub}  
\end{equation}
while our definition of $H_{+}$ implies that there exists 
$\beta>\lambda$ such that
\begin{equation}
	|Au|\geq\beta|u|
	\quad\quad 
	\forall u\in H_{+}\cap D(A).
	\label{hp:A-super}
\end{equation}

More precisely, the last inequality holds true with $\beta$ equal to
the smallest eigenvalue of $A$ greater than $\lambda$ (if any), or
with any $\beta>\lambda$ if the spectrum of $A$ is finite and
$\lambda$ is its maximum, in which case $H_{+}=\{0\}$.

Now let us choose a constant $\delta$ such that
$\lambda<\delta<\min\left\{\beta,(1+p)\lambda\right\}$, let us set
\begin{equation}
	r_{1}:=2\left(1+\frac{1}{\beta}\right)^{1/2}r_{0},
	\label{defn:r1}
\end{equation}
and let us assume that $r_{0}$ is small enough so that
\begin{equation}
	2r_{1}< R,
	\label{hp:small-r0-1}
\end{equation}
\begin{equation}
	2L(2r_{1})^{p}
	 \left(\frac{\sqrt{\lambda+1}}{\delta-\lambda}
	 +\frac{\sqrt{\beta+1}}{\beta-\delta}\right)
	 \leq \frac{1}{2}.
	\label{hp:small-r0-2}
\end{equation}

Let us consider the space 
$$\mathbb{X}:=\left\{g\in C^{0}\left([0,+\infty);D(A^{1/2})\right):
|g(t)|_{D(A^{1/2})}\leq r_{1}\quad\forall t\geq 0\right\}.$$

It is well-known that $\mathbb{X}$ is a complete metric space with respect 
to the distance
$$\mbox{dist}(g_{1},g_{2}):=\sup
\left\{|g_{1}(t)-g_{2}(t)|_{D(A^{1/2})}:t\geq 0\right\}.$$

For every $g\in \mathbb{X}$ we set
\begin{equation}
	\varphi_{g}(t):=f\left(v_{0}e^{-\lambda t}+g(t)e^{-\delta t}\right)
	\quad\quad
	\forall t\geq 0,
	\label{defn:phig}
\end{equation}
and we define $\varphi_{g,-}(t)$ and $\varphi_{g,+}(t)$ as the
projections of $\varphi_{g}(t)$ into $H_{-}$ and $H_{+}$,
respectively.  
Then we define $u_{g,-}:[0,+\infty)\to H_{-}$ and $u_{g,+}:[0,+\infty)\to
H_{+}$ as
\begin{eqnarray}
	u_{g,-}(t) & := & v_{0}e^{-\lambda t} -\int_{t}^{+\infty}
	e^{(s-t)A}\cdot\varphi_{g,-}(s)\,ds,
	\label{defn:ug-}  \\
	\noalign{\vspace{0.5ex}}
	u_{g,+}(t) & := & e^{-tA}w_{0}+\int_{0}^{t}
	e^{(s-t)A}\cdot\varphi_{g,+}(s)\,ds,
	\label{defn:ug+}  
\end{eqnarray}
and  $u_{g}(t):=u_{g,-}(t)+u_{g,+}(t)$. Finally, we set
\begin{equation}
	\overline{g}(t):=\left(u_{g}(t)-v_{0}e^{-\lambda t}\right)e^{\delta t}
	\quad\quad
	\forall t\geq 0.
	\label{defn:g-bar}
\end{equation}

We claim that the following three statements hold true, provided that
the smallness assumptions (\ref{hp:small-r0-1}) and
(\ref{hp:small-r0-2}) are satisfied.

\begin{itemize}
	\item \emph{Well-posedness of the construction}.  The functions
	$\varphi_{g}$ and $u_{g}$ are well-defined for every $g\in
	\mathbb{X}$.  Moreover, $u_{g}$ is a solution to the
	\emph{linear} equation
	\begin{equation}
		u_{g}'(t)+Au_{g}(t)=\varphi_{g}(t)= 
		f\left(v_{0}e^{-\lambda t}+g(t)e^{-\delta t}\right)
		\quad\quad
		\forall t\geq 0
		\label{eqn:u-phig}
	\end{equation}
	in the sense of Theorem~\ref{thmbibl:linear}.

	\item \emph{Closedness}.  We have that $\overline{g}\in \mathbb{X}$
	for every $g\in \mathbb{X}$.

	\item \emph{Contractivity}.  The map $\mathcal{F}:\mathbb{X}\to
	\mathbb{X}$ defined by $\mathcal{F}(g)=\overline{g}$ is a
	contraction.
\end{itemize}

If we prove the three claims, then the conclusion easily follows.
Indeed the contractivity implies that $\mathcal{F}$ has a fixed point,
namely there exists $g\in \mathbb{X}$ such that $\overline{g}=g$.  Thus from
(\ref{defn:g-bar}) we have that
\begin{equation}
	u_{g}(t)=v_{0}e^{-\lambda t}+\overline{g}(t)e^{-\delta t}=
	v_{0}e^{-\lambda t}+g(t)e^{-\delta t},
	\label{eqn:u-g}
\end{equation}
so that (\ref{th:limit}) follows from the boundedness of
$g$ and the fact that $\delta>\lambda$.  

Moreover, the projection of $u_{g}(0)$ into $H_{+}$, namely $u_{g,+}(0)$,
is $w_{0}$ because of~(\ref{defn:ug+}), hence the initial condition 
$u_{g}(0)$ is of the form $w_{1}+w_{0}$ for some $w_{1}\in H^{-}$.

Finally, exploiting~(\ref{eqn:u-g}) once again, we have that
$$f(u_{g}(t))=f\left(v_{0}e^{-\lambda t}+g(t)e^{-\delta t}\right),$$
so that saying that $u_{g}$ is a solution to (\ref{eqn:u-phig}) is
equivalent to saying that $u_{g}$ is a solution to (\ref{pbm:nl}).

\subparagraph{\textmd{\textit{Well-posedness of the construction}}}

From (\ref{hp:smallness}) and our definition of $\mathbb{X}$ we have
that 
$$\left|v_{0}e^{-\lambda t}+g(t)e^{-\delta
t}\right|_{D(A^{1/2})}\leq |v_{0}|_{D(A^{1/2})}e^{-\lambda t}+
|g(t)|_{D(A^{1/2})}e^{-\delta t}\leq
r_{0}e^{-\lambda t}+r_{1}e^{-\delta t}.$$

Since $\delta>\lambda$ and $r_{0}\leq r_{1}$, it follows that
\begin{equation}
	\left|v_{0}e^{-\lambda t}+g(t)e^{-\delta t}\right|_{D(A^{1/2})}\leq
	2r_{1}e^{-\lambda t}
	\quad\quad
	\forall t\geq 0.
	\label{est:ug}
\end{equation}

Due to the smallness assumption (\ref{hp:small-r0-1}), we have in
particular that 
$$\left|v_{0}e^{-\lambda t}+g(t)e^{-\delta
t}\right|_{D(A^{1/2})}< R \quad\quad
\forall t\geq 0,$$
which proves that $\varphi_{g}(t)$ is well-defined.

Setting $v=0$ into~(\ref{hp:f-plip}), from (\ref{hp:f0}) we obtain
that
$$|f(u)|\leq L|u|^{1+p}_{D(A^{1/2})}
\quad\quad
\forall u\in B_{R}.$$

Therefore, recalling that $\delta<(1+p)\lambda$, from
(\ref{defn:phig}) and (\ref{est:ug}) we deduce that
\begin{equation}
	|\varphi_{g}(t)|\leq 
	L(2r_{1})^{1+p}e^{-\delta t}
	\quad\quad
	\forall t\geq 0.	
	\label{est:phig}
\end{equation}

Let us examine the integrand in (\ref{defn:ug-}).  From
(\ref{hp:A-sub}) it turns out that the operator $e^{(s-t)A}$ is
bounded in $H_{-}$ with norm less than or equal to $e^{(s-t)\lambda}$.
Since the right-hand side of (\ref{est:phig}) is of course an estimate
also for $|\varphi_{g,-}(t)|$, we deduce that
\begin{equation}
	\left|e^{(s-t)A}\varphi_{g,-}(s)\right|\leq
	e^{(s-t)\lambda}\left|\varphi_{g,-}(s)\right|\leq
	e^{(s-t)\lambda}L(2r_{1})^{1+p}e^{-\delta s}
	\quad\quad
	\forall s\geq t\geq 0.
	\label{est:integrand-u-}
\end{equation}

Since $\delta>\lambda$, this proves that the integral in the
right-hand side of (\ref{defn:ug-}) converges for every $t\geq 0$,
hence $u_{g,-}(t)$ is well-defined.  

Let us examine the right-hand side of (\ref{defn:ug+}).  Now the
integration is over a bounded interval, and the operator $e^{(s-t)A}$
is a contraction because $s\leq t$.  Therefore, also $u_{g,+}(t)$ is 
well-defined.

Finally, both the regularity (\ref{th:reg-cont}) of $u_{g}$, and the
fact that it is a solution to (\ref{eqn:u-phig}), follow from
definitions (\ref{defn:ug-}) and (\ref{defn:ug+}).

\subparagraph{\textmd{\textit{Closedness}}}

To begin with, we observe that $\overline{g}:[0,+\infty)\to
D(A^{1/2})$ is a continuous map because of the regularity of $u_{g}$.
Now we claim that
\begin{equation}
	\left|u_{g,-}(t)-v_{0}e^{-\lambda t}\right|_{D(A^{1/2})}\leq
	L(2r_{1})^{1+p}
	\frac{\sqrt{\lambda+1}}{\delta-\lambda}e^{-\delta t}
	\quad\quad
	\forall t\geq 0,
	\label{est:ug-}
\end{equation}
\begin{equation}
	\left|u_{g,+}(t)\right|_{D(A^{1/2})}\leq
	L(2r_{1})^{1+p}
	\frac{\sqrt{\beta+1}}{\beta-\delta}e^{-\delta t}
	+\frac{r_{1}}{2}e^{-\delta t}
	\quad\quad
	\forall t\geq 0.
	\label{est:ug+}
\end{equation}

Plugging these estimates into (\ref{defn:g-bar}) we obtain that
\begin{eqnarray*}
	\left|\overline{g}(t)\right|_{D(A^{1/2})} & = & 
	\left|u_{g,-}(t)+u_{g,+}(t)-v_{0}e^{-\lambda t}\right|_{D(A^{1/2})}
	e^{\delta t}\\
	\noalign{\vspace{0.5ex}}
	 & \leq & 
	 \left|u_{g,-}(t)-v_{0}e^{-\lambda t}\right|_{D(A^{1/2})}e^{\delta t}+
	 |u_{g,+}(t)|_{D(A^{1/2})}e^{\delta t}  \\
	 \noalign{\vspace{0.5ex}}
	 & \leq & L(2r_{1})^{1+p}
	 \left(\frac{\sqrt{\lambda+1}}{\delta-\lambda}
	 +\frac{\sqrt{\beta+1}}{\beta-\delta}\right)
	 +\frac{r_{1}}{2}.	 
\end{eqnarray*}

The smallness assumption (\ref{hp:small-r0-2}) is equivalent to saying
that the right-hand side is less than or equal to $r_{1}$.  This
proves that $\overline{g}\in \mathbb{X}$.

Let us prove (\ref{est:ug-}).  From (\ref{est:integrand-u-}) we have that
\begin{eqnarray*}
	\left|u_{g,-}(t)-v_{0}e^{-\lambda t}\right| & \leq & 
	\int_{t}^{+\infty}\left|e^{(s-t)A}\varphi_{g,-}(s)\right|\,ds   \\
	 & \leq & L(2r_{1})^{1+p}
	 \int_{t}^{+\infty}e^{(s-t)\lambda}\cdot e^{-\delta s}\,ds\\
	 & = & L(2r_{1})^{1+p}\frac{1}{\delta-\lambda}
	 e^{-\delta t},
\end{eqnarray*} 
so that (\ref{est:ug-}) follows by simply remarking that
\begin{equation}
	|w|_{D(A^{1/2})}\leq\sqrt{\lambda+1}\cdot|w|
	\quad\quad
	\forall w\in H_{-}.
	\label{est:norm-H-}
\end{equation}

Let us prove (\ref{est:ug+}). To this end, we set 
$E(t):=|A^{1/2}u_{g,+}(t)|^{2}$. Its time-derivative is
\begin{eqnarray*}
	E'(t) & = & -2|Au_{g,+}(t)|^{2}+
	2\langle Au_{g,+}(t),\varphi_{g,+}(t)\rangle  \\
	 & \leq & -2|Au_{g,+}(t)|^{2}+
	 \frac{\beta-\delta}{\beta}|Au_{g,+}(t)|^{2}+
	 \frac{\beta}{\beta-\delta}|\varphi_{g,+}(t)|^{2} \\
	 & \leq & -\frac{\delta+\beta}{\beta}|Au_{g,+}(t)|^{2}+
	 \frac{\beta}{\beta-\delta}|\varphi_{g}(t)|^{2}.
\end{eqnarray*}

Let us estimate the first term using (\ref{hp:A-super}), and the second 
term using (\ref{est:phig}). We deduce that
$$E'(t)\leq -(\delta+\beta)E(t)+
\frac{\beta}{\beta-\delta}L^{2}(2r_{1})^{2+2p}e^{-2\delta t}.$$

Integrating this differential inequality we obtain that
$$E(t)\leq |A^{1/2}w_{0}|^{2}e^{-(\delta+\beta)t}+
\frac{\beta}{(\beta-\delta)^{2}}L^{2}(2r_{1})^{2+2p}e^{-2\delta t}
\quad\quad
\forall t\geq 0.$$

Since $\delta+\beta>2\delta$, and since $|A^{1/2}w_{0}|\leq r_{0}$ 
because of assumption (\ref{hp:smallness}), this easily implies that
$$\left|A^{1/2}u_{g,+}(t)\right|\leq r_{0}e^{-\delta t}+
\frac{\sqrt{\beta}}{\beta-\delta}L(2r_{1})^{1+p}e^{-\delta t}
\quad\quad
\forall t\geq 0,$$
so that (\ref{est:ug+}) follows by simply recalling our definition
(\ref{defn:r1}) of $r_{1}$, and the fact that
\begin{equation}
	|w|_{D(A^{1/2})}\leq
	\left(1+\frac{1}{\beta}\right)^{1/2}|A^{1/2}w|
	\quad\quad
	\forall w\in D(A^{1/2})\cap H_{+}.
	\label{est:norm-H+}
\end{equation}

\subparagraph{\textmd{\textit{Contractivity}}}

Let $g_{1}$ and $g_{2}$ be two elements of $\mathbb{X}$. From (\ref{hp:f-plip}) 
we have that
\begin{eqnarray*}
	|\varphi_{g_{1}}(t)-\varphi_{g_{2}}(t)| & = & 
	\left|f\left(v_{0}e^{-\lambda t}+g_{1}(t)e^{-\delta t}\right)-
	f\left(v_{0}e^{-\lambda t}+g_{2}(t)e^{-\delta t}\right)\right|
	  \\
	 & \leq & L\left(
	 \left|v_{0}e^{-\lambda t}+g_{1}(t)e^{-\delta t}\right|^{p}_{D(A^{1/2})}
	 +\left|v_{0}e^{-\lambda t}+g_{2}(t)e^{-\delta t}\right|^{p}_{D(A^{1/2})}
	 \right)\times \\
	  & & \times|g_{1}(t)-g_{2}(t)|_{D(A^{1/2})}e^{-\delta t}.
\end{eqnarray*}

Therefore, applying inequality (\ref{est:ug}) to $g_{1}$ and $g_{2}$,
we deduce that
\begin{equation}
	|\varphi_{g_{1}}(t)-\varphi_{g_{2}}(t)|\leq
	2L(2r_{1})^{p}\cdot
	\mbox{dist}(g_{1},g_{2})\cdot e^{-\delta t}
\quad\quad
\forall t\geq 0.
	\label{est:phig-12}
\end{equation}

The right-hand side of (\ref{est:phig-12}) is of course an estimate
also for the projections of $\varphi_{g_{1}}-\varphi_{g_{2}}$ into
$H_{-}$ and $H_{+}$.  In particular, for the component with respect
to $H_{-}$ we obtain that
\begin{eqnarray*}
	|u_{g_{1},-}(t)-u_{g_{2},-}(t)| & \leq & 
	\int_{t}^{+\infty}\left|e^{(s-t)A}(
	\varphi_{g_{1},-}(s)-\varphi_{g_{2},-}(s))\right|\,ds\\
	 & \leq & \int_{t}^{+\infty}e^{(s-t)\lambda}\left|
	\varphi_{g_{1}}(s)-\varphi_{g_{2}}(s)\right|\,ds  \\
	 & \leq & 2L(2r_{1})^{p}\cdot\mbox{dist}(g_{1},g_{2})
	 \int_{t}^{+\infty}e^{(s-t)\lambda}\cdot e^{-\delta s} \,ds \\
	 & \leq & 2L(2r_{1})^{p}\cdot\mbox{dist}(g_{1},g_{2})
	 \cdot\frac{1}{\delta-\lambda}e^{-\delta t}.
\end{eqnarray*}

Keeping (\ref{est:norm-H-}) into account, this proves that
\begin{equation}
	|u_{g_{1},-}(t)-u_{g_{2},-}(t)|_{D(A^{1/2})}\leq
	2L(2r_{1})^{p}\cdot\mbox{dist}(g_{1},g_{2})\cdot
	\frac{\sqrt{\lambda+1}}{\delta-\lambda}e^{-\delta t}.
	\label{est:ug12-}
\end{equation}

Let us consider now the component with respect to $H_{+}$.  To this
end, we set once again
$E(t):=|A^{1/2}(u_{g_{1},+}(t)-u_{g_{2},+}(t))|^{2}$.  Taking the
time-derivative, and arguing as we did in the previous paragraph, we
obtain that
$$E'(t)\leq -\frac{\delta+\beta}{\beta}
\left|A(u_{g_{1},+}(t)-u_{g_{2},+}(t))\right|^{2}+
\frac{\beta}{\beta-\delta}|\varphi_{g_{1}}(t)-\varphi_{g_{2}}(t)|^{2}.$$

Thus from (\ref{hp:A-super}) and (\ref{est:phig-12}) it follows that
$$E'(t)\leq -(\delta+\beta)E(t)+
\frac{\beta}{\beta-\delta}\cdot 4L^{2}(2r_{1})^{2p}\cdot
\mbox{dist}^{2}(g_{1},g_{2})\cdot e^{-2\delta t}.$$

Since now $E(0)=0$, integrating this differential inequality we
deduce that 
$$E(t)\leq 4L^{2}(2r_{1})^{2p}\cdot
\mbox{dist}^{2}(g_{1},g_{2})\cdot\frac{\beta}{(\beta-\delta)^{2}}
e^{-2\delta t}
\quad\quad
\forall t\geq 0.$$

Keeping (\ref{est:norm-H+}) into account, this proves that
\begin{equation}
	|u_{g_{1},+}(t)-u_{g_{2},+}(t)|_{D(A^{1/2})}\leq
	2L(2r_{1})^{p}\cdot\mbox{dist}(g_{1},g_{2})\cdot
	\frac{\sqrt{\beta+1}}{\beta-\delta}e^{-\delta t}.
	\label{est:ug12+}
\end{equation}

From (\ref{defn:g-bar}), (\ref{est:ug12-}) and (\ref{est:ug12+}) we
conclude that
\begin{eqnarray*}
	\left|\overline{g_{1}}(t)-\overline{g_{2}}(t)\right|_{D(A^{1/2})} & = & 
	\left|u_{g_{1}}(t)-u_{g_{2}}(t)\right|_{D(A^{1/2})}e^{\delta t}\\
	\noalign{\vspace{0.5ex}}
	 & \leq & \left|u_{g_{1},-}(t)-u_{g_{2},-}(t)\right|_{D(A^{1/2})}e^{\delta t}+
	 \left|u_{g_{1},+}(t)-u_{g_{2},+}(t)\right|_{D(A^{1/2})}e^{\delta t}\\
	 \noalign{\vspace{0.5ex}}
	 & \leq & 2L(2r_{1})^{p}
	 \left(\frac{\sqrt{\lambda+1}}{\delta-\lambda}
	+\frac{\sqrt{\beta+1}}{\beta-\delta}\right)\cdot
	\mbox{dist}(g_{1},g_{2})
\end{eqnarray*}
for every $t\geq 0$.  Taking the supremum over all $t\geq 0$, and
keeping into account the smallness assumption (\ref{hp:small-r0-2}),
we conclude that
$$\mbox{dist}(\overline{g_{1}},\overline{g_{2}})\leq
\frac{1}{2}\,\mbox{dist}(g_{1},g_{2}),$$
which proves that the map $\mathcal{F}:\mathbb{X}\to \mathbb{X}$ is a
contraction.\qed

\setcounter{equation}{0}
\section{Applications}\label{sec:applications}

Let $\Omega\subseteq\re^{n}$ be a bounded connected open set with
Lipschitz boundary (or any other condition which guarantees Sobolev
embeddings).  We consider semilinear parabolic equations of the form
\begin{equation}
	u_{t}-\Delta u+\psi(u)=0
	\label{pbm:n-psi}
\end{equation}
in $\Omega\times [0,+\infty)$, with homogeneous Neumann boundary
conditions, and initial datum $u(0)=u_{0}$. We also consider 
semilinear parabolic equations of the form
\begin{equation}
	u_{t}-\Delta u-\lambda u+\psi(u)=0
	\label{pbm:d-psi}
\end{equation}
in $\Omega\times [0,+\infty)$, with homogeneous Dirichlet boundary
conditions, and initial datum $u(0)=u_{0}$. In this case we 
assume that $\lambda\leq\lambda_{1}(\Omega)$, where 
$\lambda_{1}(\Omega)$ denotes the first eigenvalue of $-\Delta$ with 
homogeneous Dirichlet boundary conditions in $\Omega$.

The classical approach to these problems consists in setting
$H:=L^{2}(\Omega)$ and considering the operator $Au:=-\Delta u$ with
domain $D(A):=H^{2}(\Omega)$ in the Neumann case, or the operator
$Au:=-\Delta u-\lambda u$ with domain $D(A):=H^{2}(\Omega)\cap
H^{1}_{0}(\Omega)$ in the Dirichlet case.  In both cases $A$ is a
self-adjoint operator on $H$.  It is a coercive operator in the
subcritical Dirichlet case with $\lambda<\lambda_{1}(\Omega)$, but it
is just a nonnegative operator both in the Neumann case (where the
kernel of $A$ is the space of constant functions), and in the critical
Dirichlet case with $\lambda=\lambda_{1}(\Omega)$ (where the kernel of
$A$ is the eigenspace of $-\Delta$ relative to the eigenvalue
$\lambda_{1}(\Omega)$).

As for the nonlinear term, we assume that $\psi:\re\to\re$ is 
a function of class $C^{1}$ such that 
\begin{equation}
	\psi(0)=0,
	\label{hp:psi-0}
\end{equation}
and
\begin{equation}
	|\psi'(\sigma)|\leq C_{1}|\sigma|^{p}
	\quad\quad
	\forall\sigma\in\re,
	\label{hp:psi+}
\end{equation}
\begin{equation}
	\psi(\sigma)\sigma\geq C_{2}|\sigma|^{2+p}
	\quad\quad
	\forall\sigma\in\re,
	\label{hp:psi-}
\end{equation}
for suitable positive constants $C_{1}$, $C_{2}$, $p$.

Now equations (\ref{pbm:n-psi}) and (\ref{pbm:d-psi}) can be written 
in the abstract form (\ref{pbm:main}) provided that we set
\begin{equation}
	[f(u)](x):=-\psi(u(x))
	\quad\quad
	\forall x\in\Omega
	\label{defn:f-psi}
\end{equation}

Therefore, it is crucial to know when $f$ satisfies the assumptions of
our main abstract results. Since clearly $f(0)=0$, the key point is 
the verification of (\ref{hp:f-plip}), which we already know to imply
both (\ref{hp:f-lip}) and (\ref{hp:f-ho}) with $p=q$. The 
verification of (\ref{hp:f-plip}) is quite standard, and requires 
that $D(A^{1/2})$ is contained into $L^{2+2p}(\Omega)$, which in turn 
is equivalent to the Sobolev embedding $H^{1}(\Omega)\subseteq 
L^{2+2p}(\Omega)$.

Here we skip the details, for which the interested reader is referred 
to Section~4.1 of~\cite{GGH:sol-lentes}. The final result is the 
following.

\begin{prop}\label{prop:f}
	Let $\Omega\subseteq\re^{n}$ be a bounded connected open set with
	Lipschitz boundary.  Let $p$ be a positive exponent, with no
	further restriction if $n\in\{1,2\}$, and $p\leq 2/(n-2)$ if
	$n\geq 3$.  Let $H:=L^{2}(\Omega)$, and let $A$ be the operator
	associated to the Neumann or Dirichlet problem.  For every $R>0$,
	let $B_{R}$ be the open ball in $D(A^{1/2})$, defined as in
	Theorem~\ref{thmbibl:nl}.  Let $\psi:\re\to\re$ be a function of
	class $C^{1}$ satisfying~(\ref{hp:psi-0}) and~(\ref{hp:psi+}).
	
	Then there exist $R>0$ and $L$ such that (\ref{defn:f-psi}) 
	defines a function $f:B_{R}\to H$ satisfying (\ref{hp:f-plip}).
\end{prop}

We point out that in Proposition~\ref{prop:f} above we need no
assumption on the sign or monotonicity of $\psi$.  Sign and
monotonicity conditions are important when looking for global
solutions for all initial data $u_{0}\in L^{2}(\Omega)$.  For the
convenience of the reader, we quote the following well-known result
for the case with the ``right sign''.

\begin{thmbibl}[Right sign -- Global existence]\label{thmbibl:right}
	Let $\Omega\subseteq\re^{n}$ be a bounded connected open set with
	Lipschitz boundary.  Let $p>0$, and let $\psi:\re\to\re$ be a
	continuous nondecreasing function satisfying~(\ref{hp:psi-}).
	
	Then the following statement applies to both the Neumann problem
	for equation~(\ref{pbm:n-psi}), and to the Dirichlet problem for
	equation~(\ref{pbm:d-psi}) with $\lambda\leq\lambda_{1}(\Omega)$.
	For every $u_{0}\in L^{2}(\Omega)$, the problem admits a unique
	global solution
	\begin{equation}
		u\in C^{0}\left([0,+\infty), H\right)\cap
		C^{0}\left((0,+\infty), D(A)\right)\cap
		C^{1}\left((0,+\infty), H\right),
		\label{th:appl-reg}
	\end{equation}
	and there exists a constant $M$ such that
	\begin{equation}
		|u(t)|\leq
		\frac{M}{t^{1/p}}
		\quad\quad
		\forall t>0,
		\label{th:appl-infty}
	\end{equation}
	\begin{equation}
		|A^{1/2} u(t)|\leq 
		\frac{|u_{0}|}{t^{1/2}}
		\quad\quad
		\forall t>0.
		\label{th:appl-auq}
	\end{equation}
		
\end{thmbibl}

The proof of Theorem~\ref{thmbibl:right} is quite classical.  Indeed,
both the Neumann and the Dirichlet problem are the gradient-flow of a
strictly convex functional, hence global existence, regularity and
estimate (\ref{th:appl-auq}) follow, for example, from the
well-established theory of maximal monotone operators
(see~\cite{brezis}). As for the decay estimate~(\ref{th:appl-infty}), 
it holds true for all solutions of the abstract equation 
(\ref{pbm:main}), provided that the nonlinear term $f(u)$ satisfies a 
sign condition such as
\begin{equation}
	\langle f(u),u\rangle\leq -K|u|^{2+p}
	\quad\quad
	\forall u\in D(A^{1/2})
	\label{hp:appl-sign}
\end{equation}
for a suitable constant $K$. Indeed, setting as usual 
$y(t):=|u(t)|^{2}$, and computing its time-derivative as in 
(\ref{eqn:y'-bis}), we obtain that
$$y'(t)=-2|A^{1/2}u(t)|^{2}+2
\langle u(t),f(u(t))\rangle\leq
-2K|u(t)|^{2+p}=-2K[y(t)]^{1+p/2},$$
so that (\ref{th:appl-infty}) easily follows by integrating this
differential inequality.  In the concrete example,
(\ref{hp:appl-sign}) holds true whenever $\psi(\sigma)$
satisfies~(\ref{hp:psi-}).

We point out that Theorem~\ref{thmbibl:right} holds true for every
$n\geq 1$ and $p>0$.  Under the more restrictive assumptions of
Proposition~\ref{prop:f}, we can apply our theory as follows.

\begin{thm}[Right sign -- Decay rates]\label{thm:right-sign}
	Let $\Omega\subseteq\re^{n}$ be a bounded open set with Lipschitz
	boundary.  Let $p$ be a positive exponent, with no further
	restriction if $n\in\{1,2\}$, and $p\leq 2/(n-2)$ if $n\geq 3$.
	Let $\psi:\re\to\re$ be a nondecreasing function of class $C^{1}$
	satisfying~(\ref{hp:psi-0}) through~(\ref{hp:psi-}).
	
	Then the following three statements apply to both the Neumann
	problem for equation~(\ref{pbm:n-psi}), and to the Dirichlet
	problem for equation~(\ref{pbm:d-psi}) with
	$\lambda=\lambda_{1}(\Omega)$.
	\begin{enumerate}
		\renewcommand{\labelenumi}{(\arabic{enumi})}
		
		\item \emph{(Classification of decay rates)} 
		All non-zero solutions are either slow or fast in the sense 
		of Theorem~\ref{thm:main-alternative}.
		
		\item \emph{(Existence of slow solutions)} There exists a 
		nonempty open set $\mathcal{S}\subseteq L^{2}(\Omega)$ such 
		that all solutions with $u(0)\in S$ are slow.
		
		\item \emph{(Existence of fast solutions)} There exists 
		families of fast solutions in the sense of 
		Theorem~\ref{thm:main-exponential}.
		
	\end{enumerate}
		
\end{thm}

Of course the theory applies also to the Dirichlet problem with 
$\lambda<\lambda_{1}(\Omega)$, but in that case the operator is 
coercive and we have only fast solutions.

The proof of Theorem~\ref{thm:right-sign} is a straightforward
application of our abstract results.  We just point out that here the
set $\mathcal{S}$ of initial data originating slow solutions is
claimed to be open in $H$ and not just in $D(A^{1/2})$ as in
Theorem~\ref{thm:main-slow}.  This improvement relies on the
regularizing effect.  Indeed, it is enough to take the set
$\mathcal{S}'\subseteq D(A^{1/2})$ provided by
Theorem~\ref{thm:main-slow}, and then considering the set
$\mathcal{S}\subseteq H$ of all initial data of solutions which end up
in $\mathcal{S}'$ at time $t=1$.  This set in non-empty and open due
to the regularizing effect.

\begin{rmk}
	\begin{em}
		For the heat equation, it is possible to use
		Theorem~\ref{thm:main-alternative} in conjunction with
		additional properties of the heat flow such as smoothing
		effect in arbitrary Lebesgue-based Sobolev spaces to obtain a
		major improvement with respect to what was already known.  
		In next result we obtain a slow-exponential alternative for
		all values of $p$ and $n$.  We are confident that a relevant
		refinement of our techniques would yield in this general case
		also the remaining results, for example the existence of an
		\emph{open set} of slow solutions or the existence of fast
		solutions with arbitrary spectral decay rate.
		
	\end{em}
\end{rmk}

\begin{thm}[Arbitrary power functions -- The alternative]
	Let $\Omega\subseteq\re^{n}$ be a bounded open set with Lipschitz
	boundary.  Let $p$ be any positive exponent, and let
	$\psi(\sigma):=c|\sigma|^p\sigma$ for some $c>0$.
	
	Then both Neumann problem for equation~(\ref{pbm:n-psi}), and
	Dirichlet problem for equation~(\ref{pbm:d-psi}) with
	$\lambda=\lambda_{1}(\Omega)$ satisfy the following decay
	alternative: all non-zero solutions with bounded initial data 
	satisfy either (\ref{th:slow}) or (\ref{th:fast+}) for some 
	$\gamma>\lambda$.
\end{thm}

\paragraph{\textmd{\textit{Proof}}} 

Here global existence of solutions is well-known in the sense of
$L^\infty(\Omega)$ for initial data in the same space, and solutions
are in fact strong solutions with values in $D(A^{1/2})$ as a
consequence of smoothing effect.  In order to apply
Theorem~\ref{thm:main-alternative} we observe that
$$|f(u)|_{H}=\|f(u)\|_{L^{2}(\Omega)}\leq
\|u\|_{L^{\infty}(\Omega)}^{p-r}\cdot
\|u\|_{L^{2r+2}(\Omega)}^{1+r}.$$

Now along the trajectory, say for $t\geq 1$, we have a uniform bound
on $u$ in $L^{\infty}(\Omega)$.  It is therefore enough to choose
$r\in(0,p)$ small enough so that $D(A^{1/2})$ is contained into
$L^{2r+2}(\Omega)$ and obtain that
$$|f(u)|\leq c_{2}|u|_{D(A^{1/2})}^{1+r}.$$

Then Theorem ~\ref{thm:main-alternative} is applicable and we obtain
the alternative with $p$ replaced by $r$ in the lower estimate of slow
solutions, but from the results of either~\cite {ben-arbi} for
Neumann's case or~\cite{bah-d} for Dirichlet's case we know that any
solution which does not satisfy the ``slow condition'' corresponding to
$p$ tends to $0$ faster than any negative power of $t$.  The result
follows immediately by exclusion.\qed
\bigskip

When $\psi$ has the wrong sign, for example when
$\psi(\sigma)=-|\sigma|^{p}\sigma$, global existence can fail
(see~\cite{fujita}).  Nevertheless, exploiting the so-called potential
well, one can obtain global existence for all initial data which are
small enough with respect to the norm of $D(A^{1/2})$.  This technique
requires that the operator is coercive and controls the nonlinear term
(which means Sobolev embeddings).  Since the coerciveness of the
operator is essential, this theory applies neither to the
Neumann case, nor to the critical Dirichlet case.  In other words, the
potential well applies only to the subcritical Dirichlet case, in
which case we obtain the following result.

\begin{thm}[Wrong sign, with potential well]
	Let $\Omega\subseteq\re^{n}$ be a bounded open set with Lipschitz
	boundary.  Let $p$ be a positive exponent, with no further
	restriction if $n\in\{1,2\}$, and $p\leq 2/(n-2)$ if $n\geq 3$.
	Let $\psi:\re\to\re$ be a function of class $C^{1}$
	satisfying~(\ref{hp:psi-0}) and~(\ref{hp:psi+}).  Let us
	consider the Dirichlet problem for equation (\ref{pbm:d-psi}) with
	$\lambda<\lambda_{1}(\Omega)$.
	
	Then there exists $R>0$ such that, for every $u_{0}\in B_{R}$
	(defined as in Theorem~\ref{thmbibl:nl}), the problem has a unique
	global solution, which satisfies (\ref{th:appl-reg}),
	(\ref{th:appl-infty}), and (\ref{th:appl-auq}).
	
	Moreover, every non-zero solution is fast in the sense of
	Theorem~\ref{thm:main-alternative}, and there exist families of
	fast solutions parametrized in the sense of
	Theorem~\ref{thm:main-exponential}.
\end{thm}

When there is no potential well, Theorem~\ref{thm:main-alternative}
keeps on classifying all possible decay rates of solutions which
exist globally and decay.  On the other hand, nothing in this case
guarantees decay, or even global existence, of solutions.

Nevertheless, there is one notable exception.
Theorem~\ref{thm:main-exponential} provides families of global
solutions with exponential decay without assuming neither the
coercivity of the operator, nor sign conditions on the nonlinear term.
Therefore, even in the Neumann case and in the critical Dirichlet
case, we obtain the following existence result.  

\begin{thm}[Wrong sign, without potential well]\label{thm:wrong-sign}
	Let $\Omega\subseteq\re^{n}$ be a bounded open set with Lipschitz
	boundary.  Let $p$ be a positive exponent, with no further
	restriction if $n\in\{1,2\}$, and $p\leq 2/(n-2)$ if $n\geq 3$.
	Let $\psi:\re\to\re$ be a function of class $C^{1}$
	satisfying~(\ref{hp:psi-0}) and~(\ref{hp:psi+}).
	
	Let us consider the Neumann problem for equation~(\ref{pbm:n-psi})
	or the Dirichlet problem for equation~(\ref{pbm:d-psi}) with
	$\lambda=\lambda_{1}(\Omega)$.
	
	Then there exist families of fast solutions parametrized in the sense of
	Theorem~\ref{thm:main-exponential}.
		
\end{thm}

We conclude by pointing out that our abstract results apply also to 
equations with second order operators with non-constant coefficients, 
or with higher order operators such as $\Delta^{2}$. We also allow 
more general nonlinear terms depending on $x$ and $t$, or even 
non-local nonlinear terms.

\label{NumeroPagine}

\end{document}